\documentclass[12pt, reqno]{amsart}
\usepackage[foot]{amsaddr}
\usepackage[margin={1.25in, 1.25in}]{geometry} 
\usepackage{amssymb, amsmath, amsfonts}
\numberwithin{equation}{section} 

\usepackage{lipsum}
\makeatletter
\let\@afterindenttrue\@afterindentfalse
\@afterindentfalse
\makeatother

\usepackage{enumitem}
\setlist[enumerate]{leftmargin=*} 

\usepackage{hyperref}
\hypersetup{
	colorlinks=true,
	linkcolor=red,
	citecolor=blue,
	filecolor=cyan,      
	urlcolor=magenta
}
\usepackage{bm}

\usepackage{titlesec} 
\titleformat{\section}{\centering\scshape}{\thesection.}{0.5em}{}
\titleformat{\subsection}{\bfseries\upshape}{\bfseries\upshape\thesubsection.}{0.5em}{} 
\titleformat{\subsubsection}{\itshape}{\upshape\thesubsubsection.}{0.5em}{} 

\usepackage{tikz-cd} 

\newcommand\genfd{{\bm k}}
\newcommand\op{\mathrm{op}}
\newcommand\id{\mathrm{id}}
\newcommand\KK{\mathcal{K}}

\usepackage{mathabx}
\newcommand\btr{\blacktriangleright}
\newcommand\btl{\blacktriangleleft}

\usepackage{amsthm}
\newtheoremstyle{remark}{}{}{\upshape}{}{\bfseries}{.}{0.5em}{}

\newtheorem{theorem}{Theorem}[section]

\newtheorem{corollary}{Corollary}[section]
\newtheorem{proposition}{Proposition}[section]
\newtheorem{definition}{Definition}[section]

\theoremstyle{remark}

\newtheorem{remark}{Remark}[section]

\begin{document}
	\title{Scalar extension Hopf algebroids}
	\author{Martina Stoji\'c} 
	\address{Department of Mathematics, Faculty of Science, University of Zagreb, Bijeni\v{c}ka cesta~30, Zagreb 10000, Croatia}
	\email{stojic@math.hr} 

	\begin{abstract} \vspace{2em}
		Given a Hopf algebra $H$, Brzezi\'nski and Militaru have shown that each braided commutative Yetter--Drinfeld $H$-module algebra $A$ gives rise to an associative $A$-bialgebroid structure on the smash product algebra $A \sharp H$. They also exhibited an antipode map making $A\sharp H$ the total algebra of a Lu's Hopf algebroid over $A$. However, the published proof that the antipode is an antihomomorphism covers only a special case. In this paper, a complete proof of the antihomomorphism property is exhibited. Moreover, a new generalized version of the construction is provided. Its input is a compatible pair $A$ and $A^{\op}$ of braided commutative Yetter--Drinfeld $H$-module algebras, and output is a symmetric Hopf algebroid $A\sharp H \cong H\sharp A^{\op}$ over~$A$. This construction does not require that the antipode of $H$ is invertible. 
	\end{abstract}

	\maketitle

\section{Introduction}
\subsection{Scalar extension Lu Hopf algebroids}
\subsubsection{Hopf algebroids in the sense of Lu}
	
For a possibly noncommutative associative algebra~$A$, associative $A$-bi\-al\-ge\-bro\-ids were introduced by M.~Takeuchi in the formalism of $\times_A$-bialgebras \cite{takeuchi} and later by J-H.~Lu in the present form in~\cite{lu}. Motivated by (quantization of) groupoids in Poisson geometry, Lu focused on Hopf $A$-algebroids~$\mathcal{K}$, that is $A$-bialgebroids admitting an antipode, an antihomomorphism $\tau\colon\mathcal{K}\to\mathcal{K}$ sa\-tis\-fy\-ing some axioms (later nontrivially modified by other authors). She exhibited a class of examples, interpreted as quantum (or noncommutative) analogues of transformation groupoids. Given a finite-dimensional Hopf algebra $H$ and an algebra $A$ in the monoidal ca\-te\-gory of modules over the Drinfeld double $D(H)$, the smash product algebra $A\sharp H$ is a total algebra of a Lu's Hopf $A$-algebroid if $A$ satisfies a braided commutativity condition stated in terms of the canonical element in~$D(H)$ (R-matrix). The construction can be essentially extended to some infinite-dimensional cases, provided the $R$-matrix is still defined in a completed variant of~$D(H)$.

\subsubsection{Yetter--Drinfeld modules}
	
For a finite-dimensional Hopf algebra $H$, the ca\-te\-go\-ry of $D(H)$-modules is mo\-no\-i\-dally equivalent to the Drinfeld-Majid center~\cite{drinfeld} $\mathcal{Z}({}_H\mathcal{M})$ of the monoidal category ${}_H\mathcal{M}$ of $H$-modules and also to the monoidal category ${}_H\mathcal{Y}\mathcal{D}^H$ of Yetter--Drinfeld $H$-modules. Categories $\mathcal{Z}({}_H\mathcal{M})$ and  ${}_H\mathcal{Y}\mathcal{D}^H$ make sense for infinite-dimensional $H$ as well. The center is a purely categorical construction, while $H$-Yetter--Drinfeld mo\-du\-les are $H$-modules which are also $H$-comodules satisfying an algebraic compatibility, so-called Yetter--Drinfeld  condition~\cite{drinfeld,montg,radfordtowber,yetter}. The language of Yetter--Drinfeld mo\-du\-les is an algebraic way to extend the result of Lu beyond finite-dimensionality of $H$ and beyond the bijectivity of the antipode $S$ of $H$ used in the proofs by Lu and automatically satisfied in the finite-dimensional case.

\subsubsection{Brzezi\'nski--Militaru scalar extension} 

Brzezi\'nski and Militaru~\cite{BrzMilitaru} have proven that a smash product algebra $A\sharp H$~\cite{montg}, for a Hopf algebra $H$ and an $H$-module algebra $A$, has an $A$-bialgebroid structure provided by specific formulas if and only if $A$ is a braided commutative algebra in ${}_H\mathcal{Y}\mathcal{D}^H$. Following~\cite{bohmHbk}, 3.4.7 (where right-right Yetter--Drinfeld module algebras were used), we call such bialgebroid a scalar extension bialgebroid: it extends a Hopf algebra by a Yetter--Drinfeld module algebra. Scalar extension bialgebroids have many applications, including so called dynamization in mathematical physics~\cite{doninmudrov} and (in a formally completed version) as the underlying structure on some noncommutative phase spaces~\cite{halgoid,stojicPhD}. Article~\cite{BrzMilitaru} presents a formula and a partial proof that there is also an antipode $\tau$ in the sense of Lu on~$A\sharp H$. The axioms for the antipode $\tau$ of the bialgebroid are properly checked in~\cite{BrzMilitaru}, except for a gap in the published proof of the antihomomorphism property of $\tau$. The pro\-per\-ty has been checked only on particular binary products of algebra generators, a special case. In order to avoid a circular argument, two other nontrivial cases need a separate proof. The proof is here completed in Section~\ref{sec:ahomLuHopf}, without requiring that the Hopf algebra antipode is bijective.

\subsection{Scalar extension symmetric Hopf algebroid}
	
A symmetric Hopf $A$-al\-ge\-bro\-id involves a left $A$-bi\-al\-ge\-bro\-id and a right $A^\op$-bi\-al\-ge\-bro\-id with the same total algebra $\mathcal{K}$ and an antipode map $\tau\colon\mathcal{K}\to\mathcal{K}$ satisfying some axioms. If $H$ has a bijective antipode $S$, a recipe is given in~\cite{bohmHbk}, 4.1.5, omiting proofs, for completing a scalar extension (right) bialgebroid from~\cite{bohmHbk}, 3.4.7, to a symmetric Hopf algebroid. We exhibit appropriate variants of  Brzezi\'nski--Militaru theorem for symmetric Hopf algebroids in Section \ref{section0} without requiring bijectivity of $S$. Here a compatible pair of a left-right Yetter--Drinfeld module algebra structure on $A$ and a right-left Yetter--Drinfeld module algebra structure on $A^{\op}$ should be given as data. If the antipode~$S$ is invertible, the two Yetter--Drinfeld module algebras determine one another. We also show in Section \ref{section0} that if smash product algebras $A\sharp H$ and $H\sharp A^{\op}$ are (left and right) scalar extension bialgebroids isomorphic as algebras and with matching antipodes, then algebras $A$ and $A^{\op}$ are compatible braided commutative Yetter--Drinfeld module algebras and the two bialgebroids comprise a symmetric Hopf algebroid.
	
\subsection{Scalar extensions in other works}
	
Schauenburg~\cite{schauenburgdouble} defined a monoidal category of Yetter--Drinfeld modules over any bialgebroid. B\'alint and Szlach\'anyi~\cite{balintszlach} stu\-di\-ed the weak center of the category of modules over any symmetric Hopf algebroid $\mathcal{K}$ and proved a monoidal equivalence of the category of Yetter--Drinfeld $\mathcal{K}$-modules with the weak center of the category of $\mathcal{K}$-modules. They revisit the Brzezi\'nski--Militaru theorem in this generalized context. In particular, they show in~\cite{balintszlach} that, for a braided commutative Yetter--Drinfeld $\mathcal{K}$-module algebra, the smash product algebra $\mathcal{K}\sharp A$ is an $A$-bialgebroid and if $\mathcal{K}$ is a Frobenius Hopf algebroid, the scalar extension is (a part of the structure of) a Frobenius Hopf algebroid as well. Regarding that they resort to categorical arguments available in Frobenius case of symmetric Hopf algebroids, the antipode $\tau$ of the scalar extension is not explicitly treated in their proof. In particular, it is not shown if the scalar extension $A$-bialgebroid is (a part of) a symmetric Hopf algebroid for general $\mathcal{K}$ and~$A$.

In her thesis~\cite{stojicPhD} on internal\footnote{Internal $A$-bialgebroids, where $A$  is a monoid in some monoidal category with equalizers commuting with the tensor product, were defined in~\cite{bohmInternal}. They generalize associative $A$-bialgebroids when vector spaces are replaced by objects in such a monoidal category. An adaptation of symmetric Hopf $A$-algebroids to the internal setup is introduced in~\cite{stojicPhD} and applied to certain completions.} Hopf algebroids in the category of filtered cofiltered vector spaces, which is introduced there to formalize certain completions involving  filtered vector spaces and their duals, the author has extensively studied scalar extension symmetric Hopf algebroids when the antipode $S$ is invertible. 
	
\section{Preliminaries}
	
\subsection{Bialgebroids and Hopf algebroids}

We assume familiarity with the notion of an (associative) $A$-bialgebroid~\cite{bohmHbk,BrzMilitaru,lu}, presenting only a compact definition. Left and right~$A$-bi\-al\-ge\-bro\-ids generalize $\genfd$-bialgebras to the case where a ground field $\genfd$ is replaced by a (unital, possibly noncommutative) associative $\genfd$-algebra $A$. In this article, $\otimes_\genfd$ is denoted by $\otimes$.

\begin{definition} \label{deflb}
	Let $A$ be a $\genfd$-algebra. A \emph{left bialgebroid} $(\KK,\mu,\alpha,\beta,\Delta,\epsilon)$ over $A$, called the \emph{base algebra}, consists of  
	\begin{enumerate}
		\item a unital associative algebra $(\KK,\mu)$, called the \emph{total algebra}
		\item a homomorphism \emph{source} $\alpha\colon A\to \KK$ and an antihomomorphism \emph{target} ${\beta\colon A\to \mathcal{K}}$ whose images mutually commute (therefore defining an $A$-bimodule structure on $\mathcal{K}$ by $a.k.b = \alpha(a)\beta(b) k$ for $a,b\in A, k\in\KK$)
		\item a comonoid structure $(\KK,\Delta\colon \KK\to\KK\otimes_A\KK,\epsilon\colon\KK\to A)$ in the monoidal category of $A$-bimodules (in other words, $(\KK,\Delta,\epsilon)$ is an $A$-coring)
	\end{enumerate} 
	such that 
	\begin{enumerate} \setcounter{enumi}{3}
		\item the map $\KK\otimes(\KK\otimes\KK)\to\KK\otimes_A\KK$, $k\otimes (g\otimes h)\mapsto \Delta(k)(g\otimes h)$ factorizes through a map $\KK\otimes(\KK\otimes_A\KK)\to\KK\otimes_A\KK$ and the latter map is a left unital action
		\item the map $\KK\otimes A\to A, k\otimes a \mapsto \epsilon(k\alpha(a))$ is a left unital action.
	\end{enumerate}
\end{definition}

\begin{definition} \label{defrb}
	Let $A$ be a $\genfd$-algebra. A \emph{right bialgebroid} $(\KK,\mu,\alpha,\beta,\Delta,\epsilon)$ over $A$, called the \emph{base algebra}, consists of  
	\begin{enumerate}
		\item a unital associative algebra $(\KK,\mu)$, called the \emph{total algebra}
		\item a homomorphism \emph{source} $\alpha\colon A\to \KK$ and an antihomomorphism \emph{target} ${\beta\colon A\to \mathcal{K}}$ whose images mutually commute (therefore defining an $A$-bimodule structure on $\mathcal{K}$ by $b.k.a =  k \alpha(a)\beta(b)$ for $a,b\in A, k\in\KK$)
		\item a comonoid structure $(\KK,\Delta\colon \KK\to\KK\otimes_A\KK,\epsilon\colon\KK\to A)$ in the monoidal category of $A$-bimodules (in other words, $(\KK,\Delta,\epsilon)$ is an $A$-coring)
	\end{enumerate} 
	such that 
	\begin{enumerate} \setcounter{enumi}{3}
		\item the map $(\KK\otimes\KK)\otimes \KK\to\KK\otimes_A\KK$, $(g\otimes h)\otimes k \mapsto (g\otimes h)\Delta(k)$ factorizes through a map $(\KK\otimes_A\KK)\otimes\KK\to\KK\otimes_A\KK$ and the latter map is a right unital action
		\item the map $ A \otimes \KK\to A, a\otimes k \mapsto \epsilon(\alpha(a)k)$ is a right unital action. 
	\end{enumerate}
\end{definition}

An $A$-bialgebroid structure on an algebra $B$ which is a ring over $A\otimes A^\op$ (\cite{schauenburg}) is precisely the datum needed to construct a monoidal structure on the category of $B$-modules such that the forgetful functor from the category of $B$-modules to the monoidal category of $A$-bimodules is strict monoidal. This conceptual characterization explains why several other algebraic definitions like $\times_A$-algebras of Takeuchi~\cite{takeuchi} and a bialgebroid with anchor of Xu~\cite{xu} (unfortunately called Hopf algebroid in~\cite{xu}) are equivalent to the Lu's definition (as shown in~\cite{BrzMilitaru}).

Hopf algebroids~\cite{bohmHbk,BrzMilitaru,lu} are generalizations of Hopf algebras, roughly in the same relation to groupoids as Hopf algebras are to groups. The total algebra $\KK$ is a ge\-ne\-ra\-lization of a function algebra on the space of morphisms of a groupoid and the base algebra $A$ is a generalization of a function algebra on the space of objects of the same groupoid. Hopf algebroids are usually formalized as bialgebroids possessing a version of an antipode. More than one variant of the notion appeared, including one by Schauenburg~\cite{schauenburg} (building on Takeuchi's notion of $\times_A$-bialgebras), B\"ohm~\cite{bohmnew} (in the case of invertible antipode) and B\"ohm and Szlach\'anyi~\cite{bohmszlach} (symmetric Hopf algebroid). In the first part of the paper, Section \ref{sec:ahomLuHopf}, we are concerned with the Hopf algebroids in the sense of Lu~\cite{lu}.
	
\begin{definition}\label{def:LuHopf} 
	Let $A$ be a $\genfd$-algebra. A \emph{Lu Hopf algebroid}~\cite{lu} over $A$ is a left associative $A$-bi\-al\-ge\-bro\-id $(\KK,\mu,\alpha,\beta,\Delta,\epsilon)$ with an antipode map $\tau\colon\KK\to\KK$, which is a linear antiautomorphism satisfying 
		\begin{align*}
		\tau\circ \beta & = \alpha,
		\\
		\mu\circ (\id\otimes\tau)\circ \gamma\circ \Delta & = \alpha\circ \epsilon, 
		\\
		\mu\circ (\tau\otimes_A\id)\circ \Delta &= \beta\circ \epsilon\circ \tau,
		\end{align*}
	for some linear section $\gamma \colon \mathcal{K} \otimes_A \mathcal{K} \to \mathcal{K} \otimes \mathcal{K}$ of the projection $\pi\colon \mathcal{K}\otimes \mathcal{K}\to \mathcal{K}\otimes_A \mathcal{K}$. 
\end{definition}
Somewhat arbitrary choice of the map $\gamma$, as well as the lack of symmetry,	are disadvantages of Lu's version of a Hopf algebroid. In the second part of the paper, Section~\ref{section0}, we work with symmetric Hopf algebroids, introduced by B\"ohm and Szlach\'anyi in~\cite{bohmszlach}.
	
\begin{definition} 
	Let $A$ be a $\genfd$-algebra. Denote $L:=A$ and $R:=A^\op$. A \emph{symmetric Hopf algebroid}~\cite{bohmHbk,bohmszlach} over $A$ consists of a left $L$-bialgebroid $(\KK,\mu,\alpha_L,\beta_L,\Delta_L,\epsilon_L)$, a right $R$-bialgebroid $(\KK,\mu,\alpha_R,\beta_R,\Delta_R,\epsilon_R)$ with the same underlying total algebra $(\KK,\mu)$, and an antipode map $\tau\colon\KK\to\KK$
	which is an antihomomorphism of algebras, subject to compatibilities 
		\begin{align*} 
		\alpha_L\circ\epsilon_L\circ\beta_R = \beta_R, 
		&\quad
		\beta_L\circ\epsilon_L\circ\alpha_R = \alpha_R, 
		\\
		\alpha_R\circ\epsilon_R\circ\beta_L = \beta_L,
		&\quad
		\beta_R\circ\epsilon_R\circ\alpha_L = \alpha_L,  
		\\
		(\Delta_R\otimes_L\id)\circ\Delta_L &=
		(\id\otimes_R\Delta_L)\circ\Delta_R,   
		\\
		(\Delta_L\otimes_R\id)\circ\Delta_R &=
		(\id\otimes_L \Delta_R)\circ\Delta_L ,
		\\
		\tau\circ\beta_L =\alpha_L, &\quad
		\tau\circ\beta_R = \alpha_R ,
		\\
		\mu_{\otimes'_L}\circ(\tau\otimes\id)&\circ\Delta_L = \alpha_R\circ\epsilon_R ,
		\\
		\mu_{\otimes'_R}\circ(\id\otimes \tau)&\circ\Delta_R = \alpha_L\circ\epsilon_L.
		\end{align*} 
	For clarification on maps $\mu_{\otimes'_R}$ and $\mu_{\otimes'_L}$, the reader can consult \cite{bohmHbk}.
\end{definition}
	
	
There is also a categorically defined notion of Hopf bialgebroid by Day and  Street~\cite{daystreet}. It is shown in~\cite{bohmszlach} that it contains equivalent information to that of a symmetric Hopf algebroid. Another approach somewhat reminding Lu's Hopf algebroids, using a choice of a balancing subalgebra, is recently axiomatized in~\cite{twosha}. Recasting the scalar extension Hopf algebroids to that formalism is also treated in~\cite{twosha}.

\subsection{Yetter--Drinfeld module algebras} 

\begin{definition}\label{lrydma} 
	Let $H$ be a $\genfd$-bialgebra. A \emph{left-right Yetter--Drinfeld module algebra} over $H$ is a triple $(L,\btr,\rho)$ such that 
	\begin{enumerate} 
		\item $L$ is a $\genfd$-algebra \label{a1}
		\item $(L,\btr)$ is a left $H$-module  \label{a2}
		\item $(L,\rho \colon a\mapsto a_{[0]} \otimes a_{[1]})$ is a right $H$-comodule \label{a3}
		\item action and coaction satisfy the \emph{left-right Yetter--Drinfeld condition} \label{a4}
			\begin{equation} \label{eq:YD1}
			(h_{(1)} \btr a_{[0]}) \otimes (h_{(2)} a_{[1]}) = (h_{(2)} \btr a)_{[0]} \otimes  (h_{(2)} \blacktriangleright a)_{[1]} h_{(1)}, \quad \forall h \in H, a\in L
			\end{equation}	
		\item action $\btr$ is \emph{Hopf}, that is \label{a5}
			\begin{equation}
			h\blacktriangleright (ab) = (h_{(1)}\blacktriangleright a)(h_{(2)}\blacktriangleright b), \quad h\blacktriangleright 1_L = \epsilon(h) 1_L, \quad \forall h\in H, a,b\in L
			\end{equation}
		\item coaction $\rho$ satisfies (see~\cite{BrzMilitaru}) \label{a6}
			\begin{equation}\label{eq:ma}
			(ab)_{[0]} \otimes (ab)_{[1]} = a_{[0]}b_{[0]} \otimes b_{[1]} a _{[1]}, \quad \rho(1_L) = 1_L \otimes 1_H, \quad \forall a,b\in L,
			\end{equation} 
		called the \emph{comodule algebra property} (over $H^\op$).
	\end{enumerate}
	We say that $L$ is \emph{braided commutative} if
	\begin{equation} \label{eq:bc}
		x_{[0]}(x_{[1]}\blacktriangleright a) = a x, \quad \forall a,x\in L.
	\end{equation}
\end{definition}

\begin{definition}\label{rlydma} 
	Let $H$ be a $\genfd$-bialgebra. A \emph{right-left Yetter--Drinfeld module algebra} over $H$ is a triple $(R,\btl,\lambda)$ such that 
	\begin{enumerate} 
		\item $R$ is a $\genfd$-algebra \label{b1}
		\item $(R,\btl)$ is a right $H$-module  \label{b2}
		\item $(R,\lambda\colon a\mapsto a_{[-1]}\otimes a_{[0]})$ is a left $H$-comodule	\label{b3}		
		\item action and coaction satisfy the \emph{right-left Yetter--Drinfeld condition} \label{b4}
			\begin{equation} \label{eq:YDrl1}
			a_{[-1]} h_{(1)}\otimes (a_{[0]} \btl h_{(2)}) = h_{(2)}(a \btl h_{(1)} )_{[-1]} \otimes  (a \btl h_{(1)})_{[0]}, \quad \forall h \in H, a\in R
			\end{equation}	
		\item action $\btl$ is \emph{Hopf}, that is \label{b5}
			\begin{equation}
			(ab) \btl h = 	(a \btl h_{(1)})(b \btl h_{(2)}), \quad 1_R \btl h = \epsilon(h) 1_R, \quad \forall h\in H, a,b\in R
			\end{equation}		
		\item coaction $\lambda$ satisfies \label{b6}
			\begin{equation} \label{eq:marl}
			(ab)_{[-1]} \otimes (ab)_{[0]} = b_{[-1]}a_{[-1]} \otimes a_{[0]} b_{[0]}, \quad \lambda(1_R) = 1_H\otimes 1_R, \quad \forall a,b\in R
			\end{equation} 
		called the \emph{comodule algebra property} (over $H^\op$).
	\end{enumerate}
	We say that $R$ is \emph{braided commutative} if
		\begin{equation} \label{eq:bcrl}
		(a \btl y_{[-1]})y_{[0]} = ya, \quad \forall a,y\in R.
		\end{equation} 
\end{definition}

Axioms (\ref{a2}), (\ref{a3}) and (\ref{a4}) from Definition \ref{lrydma} together comprise the axioms for $(L,\btr,\rho)$ being a \emph{left-right Yetter--Drinfeld module} over $H$, and likewise, from Definition \ref{rlydma} for $(R,\btl,\lambda)$ being a \emph{right-left Yetter--Drinfeld module} over $H$. Morphisms of Yetter--Drinfeld $H$-modules are $H$-module morphisms which are also $H$-comodule morphisms. Left-right Yetter--Drinfeld $H$-modules form a category ${}_H\mathcal{Y}\mathcal{D}^H$ carrying canonical monoidal structure with a pre-braiding \cite{radfordtowber}, which is invertible (that is, a braiding) if $H$ is a Hopf algebra. Likewise, there is a pre-braided monoidal cate\-gory ${}^H\mathcal{Y}\mathcal{D}_H$ of right-left Yetter--Drinfeld $H$-modules. 
If $H$ is a Hopf algebra with an invertible antipode, these two categories are equivalent. In more detail, if $H$ is a Hopf algebra with an antipode $S$, a left Hopf $H$-action~$\blacktriangleright$ on algebra $A$ induces a right Hopf $H$-action~$\blacktriangleleft$ on $A^\op$ by $ a\blacktriangleleft h = S(h)\blacktriangleright a$. Similarly, right $H$-coactions on $A$ with comodule algebra property can be induced from left ones on $A^\op$ by $a_{[0]}\otimes a_{[1]} = a_{[0]} \otimes S(a_{[-1]})$. If $S$ is invertible, the same identities also define left actions from right actions and left coactions from right coactions. These cor\-res\-pon\-den\-ces 
 induce equivalence of monoidal categories ${}^H\mathcal{Y}\mathcal{D}_H\cong {}_H\mathcal{Y}\mathcal{D}^H$ for any Hopf algebra $H$ with an invertible antipode.

Axioms (\ref{a1}), (\ref{a2}) and (\ref{a5}) from Definition \ref{lrydma} together say that $(L,\btr)$ is a \emph{left module algebra} over $H$. Given any left module algebra $(L,\btr)$ over $H$, the vector space $L\otimes H$ carries a structure of an associative $\genfd$-algebra with multiplication bilinearly extending formula
	\begin{equation}
	(a\otimes h)\cdot( a'\otimes h') = \sum a (h_{(1)}\blacktriangleright a')\otimes h_{(2)} h', \quad  h,h'\in H, a,a' \in L
	\end{equation} 
and unit $1_L\otimes 1_H$. This algebra is called the \emph{smash product algebra} (\cite{montg}) and is denoted $L\sharp H$. It comes along with the canonical algebra monomorphisms $L\cong L\otimes\genfd\hookrightarrow L\sharp H$ and $H\cong \genfd\otimes H\hookrightarrow L\sharp H$. The images of these two embeddings are denoted $L\sharp 1$ and $1\sharp H$. A general element $a\otimes h$ of $L\sharp H$ is denoted by $a\sharp h$. Likewise, axioms (\ref{b1}), (\ref{b2}) and (\ref{b5}) from Definition \ref{rlydma} together say that $(R,\btl)$ is a \emph{right module algebra} over $H$. Similarly, any right module algebra $(R,\btl)$ over $H$ induces the \emph{smash product algebra	$H\sharp R$} with multiplication bilinearly extending formula
	\begin{equation}
	(h\otimes a)\cdot(h'\otimes a') = \sum h h'_{(1)}\otimes (a\blacktriangleleft h'_{(2)}) a', \quad  h,h'\in H, a,a' \in R
	\end{equation} 
and unit $1_H\otimes 1_R$.
 
Written in terms of the corresponding smash product algebras $L\sharp H$ and $H\sharp R$, respectively, the left-right and the right-left Yetter--Drinfeld condition is
	\begin{align}
	h \cdot \rho(a) 
	&= \rho(h_{(2)} \btr a) \cdot  h_{(1)}, \quad  \forall a\in L, h\in H, \label{eq:YD} 
	\\
	\lambda(a) \cdot h
	&= h_{(2)} \cdot \lambda(a \btl h_{(1)} ), \quad  \forall a\in R, h\in H. \label{eq:YDrl}
	\end{align}

Axioms (\ref{a1}), (\ref{a3}) and (\ref{a4}) from Definition \ref{lrydma} and likewise from Definition \ref{rlydma} respectively define a \emph{right comodule algebra} $(L,\rho)$ over $H^\op$ and a \emph{left comodule algebra} $(R,\lambda)$ over~$H^\op$.

Therefore, a left-right Yetter--Drinfeld $H$-module algebra is a left-right Yetter--Drin\-feld $H$-module and an algebra that is a left $H$-module algebra and a right $H^\op$-comodule algebra, and likewise for the right-left case.

\subsection{Brzezi\'nski--Militaru scalar extension Hopf algebroid}

If $A$ is a braided commutative left-right Yetter--Drinfeld module algebra over a Hopf algebra $H$, then the smash product $\mathcal K = A\sharp H$ is a total algebra of a 	Hopf $A$-algebroid, called a {scalar extension Hopf algebroid}. For a Lu's Hopf algebroid this is proven in~\cite{BrzMilitaru}, modifying slightly an earlier construction of Lu~\cite{lu}, Section~5, where, instead of Yetter--Drinfeld modules, closely related modules over Drinfeld double $D(H)$ are considered. In the proof of the theorem in~\cite{BrzMilitaru}, the antihomomorphism property of the Hopf algebroid antipode $\tau$ is checked on binary products of generators of the form $a\sharp 1 \cdot 1\sharp h = a\sharp h$, and in this case it holds by the definition of $\tau$. The complete proof of the property, which by linearity amounts to showing that
	\begin{equation*}
	\tau(a\sharp h \cdot b\sharp g) = \tau(b\sharp g) \tau(a\sharp h)
	\end{equation*} 
for general elements of the form $a\sharp h, b\sharp g \in A\sharp H$, is omitted. Antihomomorphism pro\-per\-ty is checked in author's dissertation~\cite{stojicPhD}, but with the additional assumption of bijectivity of the Hopf algebra antipode $S$. We give the proof of the antihomomorphism property in Section \ref{sec:ahomLuHopf}, without additional assumptions, hence we complete the proof of the theorem in \cite{BrzMilitaru}.
	
In the theorem in \cite{BrzMilitaru}, the $A$-bimodule structure of $A\sharp H$ is determined by the 	source and target maps
	\begin{equation}
	\alpha(a) = a\sharp 1,\quad \beta(a) = a_{[0]}\sharp a_{[1]}, \quad \text{ for } a \in A,
	\end{equation}
which commute by the braided commutativity of $A$. The comonoid structure of $A\sharp H$ is given by 
	\begin{equation}
	\Delta_{A\sharp H}(a\sharp h) = (a\sharp h_{(1)}) \otimes_A (1\sharp h_{(2)}),
	\quad
	\epsilon_{A\sharp H}(a\sharp h) = a\epsilon_H(h), \quad \text{ for } a \in A \text{ and } h\in H. 
	\end{equation}
Finally, the antipode $\tau$ for the Lu's Hopf algebroid is (cf.~\cite{BrzMilitaru,lu})  
	\begin{equation}\label{eq:tause}
	\tau(a\sharp h) = S(h) S^2(a_{[1]}) \cdot a_{[0]},  \quad \text{ for } a \in A \text{ and } h\in H,
	\end{equation}
where the dot denotes the multiplication on $A\sharp H$.
	
\section{Brzezi\'nski--Militaru scalar extension antipode}	\label{sec:ahomLuHopf}
\subsection{Proof of the antihomomorphism property}
	
Here we prove that the antipode defined in the Brzezi\'nski--Militaru theorem in \cite{BrzMilitaru} is an antihomomorphism.
	
\begin{theorem} \label{BM} 
	Let $H$ be a Hopf algebra with an antipode $S$ and let $(A,\btr,\rho)$ be a braided commutative left-right Yetter--Drinfeld $H$-module algebra. Then the linear map $\tau \colon A\sharp H \to A\sharp H$ defined by 
		\begin{equation*} 
		\tau(a\sharp h) = S(h) S^2(a_{[1]}) \cdot a_{[0]}
		\end{equation*} 
	is an antihomomorphism.
\end{theorem}

\begin{proof}
	After we prove the statement for special cases 
	\begin{enumerate}
		\item $\tau( h \cdot g) = \tau(g) \tau(h), \qquad \forall h,g\in H$ \label{11}
		\item $\tau(a\sharp h)  = \tau(h) \tau(a), \qquad \forall a \in A, h\in H$ \label{12}
		\item $\tau(1\sharp h \cdot a\sharp 1) = \tau(a) \tau(h), \qquad \forall a\in A, h\in H$ \label{13}
		\item $\tau(a \cdot b) = \tau(b) \tau(a), \qquad \forall a,b\in A$, \label{14}
	\end{enumerate}
	the general statement follows directly: 
		\begin{align*}
		\tau(a\sharp h \cdot b\sharp g ) & = \tau(a (h_{(1)}\blacktriangleright b) \sharp h_{(2)} g) & \\
		& =  \tau(h_{(2)}g) \tau(a (h_{(1)}\blacktriangleright b)) & \text{by (\ref{12})}  \\
		& =\tau(g)\tau(h_{(2)}) \tau(h_{(1)}\blacktriangleright b) \tau(a) & \text{by (\ref{11}), (\ref{14})} \\ 
		& =\tau(g) \tau((h_{(1)}\blacktriangleright b)\sharp h_{(2)}) \tau(a) & \text{by (\ref{12})} \\ 
		& =\tau(g) \tau(1\sharp h \cdot b \sharp 1) \tau(a) & \\ 
		& =\tau(g) \tau ( b ) \tau( h )  \tau(a) & \text{by (\ref{13})} \\ 
		& =\tau(b\sharp g) \tau(a \sharp h ) & \text{by (\ref{12}).}
		\end{align*}
		
	Equation (\ref{11})  holds because the antipode $S$ is an antihomomorphism and the equation (\ref{12}) follows directly from the definition of $\tau$. 
				
	We now prove (\ref{13}). By the definition of $\tau$, we have
		\begin{align*}
		\tau(1\sharp h \cdot a\sharp 1) & = \tau((h_{(1)} \blacktriangleright a) \sharp h_{(2)}) & \\
		& = S(h_{(2)}) S^2((h_{(1)} \blacktriangleright a)_{[1]}) \cdot (h_{(1)} \blacktriangleright a)_{[0]}, & \\
		\tau(a) \cdot \tau(h) & = S^2(a_{[1]}) \cdot a_{[0]} \cdot S(h). &
		\end{align*}
	We have to prove the equality of these two expressions. First, by the left-right Yetter--Drinfeld property (\ref{eq:YD}) for $h \in H$ and $a \in A$, and by switching $x\otimes y \mapsto y \otimes x$, we have
		\begin{equation*}
		(h_{(2)} \blacktriangleright a)_{[1]} h_{(1)} \otimes (h_{(2)} \blacktriangleright a)_{[0]} = h_{(2)} a_{[1]} \otimes( h_{(1)} \blacktriangleright a_{[0]})
		\end{equation*}
	and from this we get
		\begin{equation*}
		h_{(4)} \otimes (h_{(3)} \blacktriangleright a)_{[1]} h_{(2)} \otimes h_{(1)} \otimes (h_{(3)} \blacktriangleright a)_{[0]} = h_{(4)} \otimes h_{(3)} a_{[1]} \otimes h_{(1)} \otimes( h_{(2)} \blacktriangleright a_{[0]})
		\end{equation*}
	by using it for $h_{(2)}$ and $a$ in the second and fourth component of the tensor product with other components the same. We now apply the map $S \otimes S^2 \otimes S\otimes \mathrm{id} $ to the tensor product and multiply all components. On the left-hand side, we get 
		\begin{equation*}
		S(h_{(4)}) \cdot S^2((h_{(3)} \blacktriangleright a)_{[1]}) S^2( h_{(2)}) \cdot S(h_{(1)}) \cdot (h_{(3)} \blacktriangleright a)_{[0]}, 
		\end{equation*}
	which is equal to the left-hand side of the equation we are set to prove,
		\begin{equation*}
		S(h_{(2)}) \cdot S^2((h_{(1)} \blacktriangleright a)_{[1]})  \cdot (h_{(1)} \blacktriangleright a)_{[0]} = \tau(1\sharp h   \cdot a \sharp 1).
		\end{equation*}
	On the right-hand side, we get
		\begin{equation*}
		S(h_{(4)}) \cdot S^2( h_{(3)}) S^2( a_{[1]} ) \cdot S(h_{(1)}) \cdot ( h_{(2)} \blacktriangleright a_{[0]}),
		\end{equation*}
	which equals
		\begin{align*}
		S^2( a_{[1]} ) \cdot S(h_{(1)}) \cdot ( h_{(2)} \blacktriangleright a_{[0]}) & =  S^2( a_{[1]} ) \cdot (S(h_{(1)})_{(1)}  h_{(2)} \blacktriangleright a_{[0]}) S(h_{(1)})_{(2)} \\
		& = S^2( a_{[1]} ) \cdot (S(h_{(1)(2)})  h_{(2)} \blacktriangleright a_{[0]}) S(h_{(1)(1)}) \\
		& = S^2( a_{[1]} ) \cdot (S(h_{(2)})  h_{(3)} \blacktriangleright a_{[0]}) S(h_{(1)}) \\
		& = S^2( a_{[1]} ) \cdot a_{[0]} \cdot S(h) \\
		& = \tau(a) \tau(h).
		\end{align*}
	Equation (\ref{13}) is proven.
		
	To prove (\ref{14}), first we prove the auxiliary statement  
		\begin{equation}\label{aux}
		\tau(a) \cdot c\sharp 1 = c\sharp 1 \cdot \tau(a), \qquad \forall a,c\in A,
		\end{equation}
	that is, that the elements in $\tau(A\sharp 1)$ commute with the elements in $A\sharp 1$. We actually use the fact that, in the left bialgebroid $A\sharp H$, elements of the image of the target map $\beta$ commute with elements of the image of the source map $\alpha$, which follows from the braided commutativity (\ref{eq:bc}),
		\begin{equation}\label{com}
		a_{[0]} \sharp a_{[1]} \cdot d = a_{[0]}  (a_{[1]} \blacktriangleright d ) \sharp a_{[2]} = d \cdot a_{[0]} \sharp a_{[1]} , \qquad \forall a,d \in A.
		\end{equation}
	We have
		\begin{align*}
		\tau(a) \cdot c & =  S^2(a_{[1]}) \cdot a_{[0]} \cdot c &  \\
		& = S^2(a_{[3]}) \cdot a_{[0]} \sharp a_{[1]} S(a_{[2]}) \cdot c &  \\
		& = S^2(a_{[3]}) \cdot a_{[0]} \sharp a_{[1]} \cdot  (S(a_{[2]})_{(1)} \blacktriangleright c) \sharp S(a_{[2]})_{(2)} &  \\
		& = S^2(a_{[3]}) \cdot   (S(a_{[2]})_{(1)} \blacktriangleright c) \cdot a_{[0]} \sharp a_{[1]} \cdot S(a_{[2]})_{(2)} & \text{by (\ref{com})}\\
		& = (S^2(a_{[3]})_{(1)}  S(a_{[2](2)}) \blacktriangleright c) \sharp S^2(a_{[3]})_{(2)} \cdot a_{[0]} \sharp a_{[1]} S(a_{[2](1)}) & \\
		& = (S^2(a_{[3](1)})  S(a_{[2](2)}) \blacktriangleright c) \sharp S^2(a_{[3](2)}) \cdot a_{[0]} \sharp a_{[1]} S(a_{[2](1)}) & \\
		& = (S^2(a_{[4]})  S(a_{[3]}) \blacktriangleright c) \sharp S^2(a_{[5]}) \cdot a_{[0]} \sharp a_{[1]} S(a_{[2]}) & \\
		& = c \cdot S^2(a_{[1]}) \cdot a_{[0]} & \\
		& = c \cdot \tau(a).
		\end{align*}  
	By this, auxiliary statement (\ref{aux}) is proven. Now the statement (\ref{14}) follows easily since 
		\begin{align*}
		\tau(b) \cdot \tau(a) & =  S^2(b_{[1]}) \cdot b_{[0]} \cdot \tau(a) &  \\
		& = S^2(b_{[1]}) \cdot \tau(a) \cdot b_{[0]} & \text{by (\ref{aux})} \\
		&= S^2(b_{[1]}) S^2(a_{[1]})  \cdot a_{[0]} b_{[0]}
		\end{align*}  
	and 
		\begin{align*}
		\tau(ab)  & =  S^2((ab)_{[1]}) \cdot (ab)_{[0]}  &  \\
		& = S^2(b_{[1]}a_{[1]}) \cdot a_{[0]} b_{[0]} & \text{by  (\ref{eq:ma})} \\
		& = S^2(b_{[1]}) S^2(a_{[1]}) \cdot a_{[0]} b_{[0]} .&\\
		\end{align*}  
	The proof is complete.
\end{proof}

\section{Scalar extension symmetric Hopf algebroid}\label{section0}

\subsection{Overview of results} Any Lu Hopf algebroid is also a symmetric Hopf algebroid if the Hopf algebra antipode $S$ is bijective \cite{bohmHbk}, therefore, in that case, the Brzezi\'nski--Militaru theorem produces also a symmetric Hopf algebroid. In the following subsections, we prove the similar theorem for symmetric Hopf algebroids without the assumption of bijectivity of $S$. In that case, two braided commutative Yetter--Drinfeld module algebras that are compatible in a certain way are needed to produce a scalar extension symmetric Hopf algebroid.

In Subsection \ref{sectionA}, we prove that having two braided commutative Yetter--Drinfeld module algebras $(L,\btr, \rho)$ and $(R,\btl, \lambda)$ which are compatible in a certain way, defined in Definition \ref{compphi}, the corresponding scalar extension bialgebroids are isomorphic as algebras, by a certain isomorphism $\Psi \colon L\sharp H \cong H\sharp R$ that agrees with corresponding antipodes, and they carry a canonical structure of a symmetric Hopf algebroid. This is proven in Theorem \ref{Phi} and Theorem \ref{HAHA}.
	
In Subsection \ref{sectionB}, we state analogous results in Theorem \ref{barPhi} for a similar com\-pa\-ti\-bi\-li\-ty, defined in Definition \ref{comptheta}. We also explain how both compatibilities naturally arise when considering isomorphic scalar extension bialgebroids. We then prove that these two compatibilites of Yetter--Drinfeld module algebras are equivalent and that the resulting scalar extension symmetric Hopf algebroids are the same. 

In Subsection \ref{sectionC}, with the assumption of having an isomorphism between smash product algebras $\Psi \colon L\sharp H \to H\sharp R$ that fixes $H$, identifies images of antihomomorphism coactions $\lambda, \rho$ with $L, R$ respectively, and identifies corresponding antipode maps, we prove that $(L,\btr, \rho)$ and $(R,\btl, \lambda)$ are braided commutative Yetter--Drinfeld module algebras and are compatible in both aforementioned  ways. This is proven in Theorem \ref{obrat}, which can be seen as the  theorem converse to the previous theorems.

In Subsection \ref{sectionD}, for convenience, we present all formulas for the case when the antipode $S$ is invertible.

\subsection{Compatible Yetter--Drinfeld module algebras. Scalar extension sym\-me\-tric Hopf algebroid} \label{sectionA}
	
First we present the known result about left scalar extension bialgebroids, with addi\-ti\-o\-nal version for right scalar extension bialgebroids. 
	
\begin{proposition} \label{bial}
	Let $H$ be a Hopf algebra. 
	\begin{enumerate}
		\item Let $(L,\btr)$ be a left $H$-module algebra and $(L, \rho)$ a right $H$-comodule. Then $(L,\btr,\rho)$ is a braided commutative Yetter--Drinfeld module algebra if and only if $L\sharp H$ is a left bialgebroid given structure maps \label{21}
			\begin{align*} 
			\alpha_L \colon & L \to L\sharp H &  && \alpha_L(x) &= x \sharp 1_H \\
			\beta_L  \colon & L \to L \sharp H & && \beta_L(x) &= \rho(x) =  x_{[0]} \sharp  x_{[1]}\\
			\Delta_L  \colon & L \sharp H \to L\sharp H \otimes_L L \sharp H & && \Delta_L(x \sharp f) &= x \sharp f_{(1)} \otimes_L 1_L \sharp f_{(2)} \\
			\epsilon_L \colon & L  \sharp H \to L & && \epsilon_L(x\sharp f ) &= \epsilon(f) x .
			\end{align*}
		\item Let $(R,\btl)$ be a right $H$-module algebra and $(R,\lambda)$ a left $H$-comodule. Then $(R,\btl,\lambda)$ is a braided commutative Yetter--Drinfeld module algebra if and only if $H\sharp R$ is a right bialgebroid given structure maps \label{22}
			\begin{align*} 
			\alpha_R'  \colon& R \to H\sharp R &&&  \alpha_R'(y) &= 1_H \sharp y \\
			\beta_R'  \colon& R \to H\sharp R &&& \beta_R'(y) &= \lambda(y) =  y_{[-1]} \sharp  y_{[0]}\\
			\Delta_R'  \colon & H \sharp R \to H\sharp R \otimes_R H \sharp R  &&& \Delta_R'( f \sharp y) &= f_{(1)}\sharp 1_R \otimes_R f_{(2)}\sharp y \\
			\epsilon_R' \colon& H \sharp R \to R &&& \epsilon_R'(f\sharp y) &= \epsilon(f) y.
			\end{align*}
	\end{enumerate}
\end{proposition}

\begin{proof} 
Bialgebroid definitions for parts (\ref{21}) and (\ref{22}) respectively are given in Definitions \ref{deflb} and \ref{defrb}. Part (\ref{21}) of this proposition is proven in \cite{BrzMilitaru} and part (\ref{22}) is analogous to it. We present a short overview of the proof. The antihomomorphism property of $\beta$ is equivalent to the braided commutativity and the comodule algebra property together. Commuting of elements in $\operatorname{Im}\alpha$ with elements in $\operatorname{Im}\beta$ follows from braided commutativity, therefore, the bimodule structure on the smash product algebra is obtained. Comodule structure on the smash product algebra in the category of bimodules is exhibited by using the Hopf algebra structure on $H$. Having all this, multiplication of elements in the image of $\Delta$ is well defined, that is, part (4) of the bialgebroid definition holds, if and only if Yetter--Drinfeld condition is satisfied. Part (5) of the bialgebroid definition follows from the Hopf algebra structure on~$H$. 
\end{proof}
	
Now we define a certain compatibility of braided commutative Yetter--Drinfeld module algebras and prove that each such compatible pair $L$ and $R$ gives rise to an isomorphism between the corresponding smash product algebras~$L\sharp H $ and $H\sharp R$. We then use this isomorphism to connect the scalar extension bialgebroid structures on $L\sharp H $ and $H \sharp R$ and obtain the antihomomorphism~$\tau$, which is to be a symmetric Hopf algebroid antipode.
	 
\begin{definition}\label{compphi}
Let $H$ be a Hopf algebra with an antipode $S$. Let $\phi \colon L \to R$ be an antiisomorphism of algebras. We say that a left-right braided commutative Yetter--Drinfeld $H$-module algebra $(L,\btr,\rho)$ and a right-left braided commutative Yetter--Drin\-feld $H$-module algebra $(R,\btl,\lambda)$ are \emph{compatible via} $\phi$ if their structure maps satisfy, for all $x \in L$ and $f \in H$, 
	\begin{align*}
	\phi(f \btr x) &= \phi(x) \btl Sf, \\ 
	\phi(x)_{[-1]} \otimes \phi(x)_{[0]} &= S(x_{[1]}) \otimes \phi( x_{[0]}) .
	\end{align*} 
\end{definition}

\begin{remark} 
	If Hopf algebra antipode $S$ is bijective, for a left-right Yetter--Drinfeld module algebra $L$ there always exists a unique right-left Yetter--Drinfeld module algebra $R$ compatible via $\phi$ to it, since these formulas can be used to define one structure using the other, and likewise vice versa. This is the content of Proposition~\ref{paired}, but for now, we work with compatible Yetter--Drinfeld module algebras without the assumption that $S$ is bijective.  
\end{remark}
	
\begin{remark}
	This definition of compatibility is not symmetric, i.e.\ it is possible to analogously define another compatibility via antiisomorphisms $R \to L$. Such other compatibility is given in Definition \ref{comptheta} and analogous results regarding it are presented in Theorem \ref{barPhi} in the following subsection. 
	It should be noted that the com\-pa\-ti\-bi\-li\-ty via~$\phi \colon L \to R$ is not equivalent to this other com\-pa\-ti\-bi\-li\-ty via ${\phi^{-1} \colon R \to L}$. It is actually equivalent to the other com\-pa\-ti\-bi\-li\-ty via certain antiisomorphism ${\theta \colon R \to L}$, which is defined in Theorem \ref{compa}.
\end{remark}

\begin{theorem} \label{Phi}
	Let $H$ be a Hopf algebra with an antipode $S$. Let $\phi \colon L \to R$ be an antiisomorphism of  algebras. Let $(L, \btr, \rho)$ and $(R,\btl, \lambda)$ be a left-right and a right-left braided commutative Yetter--Drinfeld $H$-module algebras  compatible via~$\phi$, that is, for all $x \in L$ and $f \in H$, 
		\begin{align}
		\phi(f \btr x) &= \phi(x) \btl Sf, \label{a}
		\\
		\phi(x)_{[-1]} \otimes \phi(x)_{[0]} &= S(x_{[1]}) \otimes \phi( x_{[0]}). \label{b} 
		\end{align} 
	Then the following holds.
	\begin{enumerate}
		\item The map $\Phi \colon L\sharp H \to H \sharp R$ defined by \label{31}
			\begin{equation}
				\Phi(x\sharp f) = Sf \sharp \phi(x)
			\end{equation}
		is an antihomomorphism of algebras.
		\item The map $\Psi \colon L\sharp H \to H\sharp R$ defined by \label{32}
			\begin{equation}
				\Psi ( x \sharp f ) = \lambda(\phi(x)) \cdot f = \phi(x)_{[-1]} \sharp \phi(x)_{[0]} \cdot f 
			\end{equation} 
		is an isomorphism of algebras with inverse $\Psi^{-1}\colon H\sharp R \to L\sharp H$,
			\begin{equation}
				\Psi^{-1}(f \sharp y) =  f S(y_{[-1]}) \cdot \phi^{-1}( y_{[0]})
			.
			\end{equation}
		Furthermore, $\Psi(L\sharp 1) = \operatorname{Im}\lambda$ and $\Psi^{-1}(1\sharp R) = \operatorname{Im}\rho$.
		\item The maps $\tau \colon L \sharp H \to L \sharp H$ and $\tau' \colon H\sharp R \to H\sharp R$ defined by 
			\begin{align}
			\tau ( x \sharp f) & = Sf S^2(x_{[1]}) \cdot x_{[0]} ,\\
			\tau' ( f \sharp y) & =   y_{[0]} \cdot S^2(y_{[-1]})Sf  
			\end{align} 
		satisfy $\Psi \circ \tau = \tau' \circ \Psi = \Phi$ and are therefore antihomomorphisms of algebras. \label{33}
	\end{enumerate}
	The maps $\Phi$, $\tau$ and $\tau'$ are bijective if and only if $S$ is bijective.
\end{theorem}
	
\begin{proof} 
	We prove (\ref{31}). Since $\phi$ and $S$ are antihomomorphisms, we have
		\begin{equation*}
		\begin{array}{rl}
			\Phi(x \sharp f \cdot d \sharp g) 
			& = \Phi((x (f_{(1)} \btr d) \sharp f_{(2)} g) \\
			& = S(f_{(2)}g) \sharp \phi(x (f_{(1)} \btr d)) \\
			& = S(g)S( f_{(2)}) \sharp \phi (f_{(1)} \btr d) \phi(x) \\
			& = S(g)S( f_{(2)}) \sharp (\phi(d) \btl S(f_{(1)})) \phi(x) \\
			& = S(g)\sharp \phi(d) \cdot S( f)\sharp \phi(x) \\
			& = \Phi(d\sharp g) \cdot \Phi(x\sharp f).
		\end{array} 
		\end{equation*}
		
	We prove (\ref{32}). Similarly as in the proof of Theorem \ref{BM}, after we prove the statement for special cases, for $x,d \in L$, $f,g \in H$, 
	\begin{enumerate}[label={(\alph*)}, ref=\alph*]
			\item $\Psi( f \cdot g) = \Psi(f) \Psi(f')$ \label{41}
			\item $\Psi(x\sharp f)  = \Psi(x) \Psi(f)$   \label{42}
			\item $\Psi(1\sharp f \cdot x\sharp 1) = \Psi(f) \Psi(x)$ \label{43}
			\item $\Psi(x \cdot d) = \Psi(x) \Psi(d)$, \label{44}
	\end{enumerate}
	the antihomomorphism property follows directly.
	Parts (\ref{41}) and (\ref{42}) follow directly from the definition. Part (\ref{44}) holds since $\phi$ and $\lambda$ are antihomomorphisms, the latter being true because $R$ is braided commutative and has the comodule algebra property:
		$$
		\begin{array}{rl}
			e_{[-1]} \sharp  e_{[0]} \cdot y_{[-1]} \sharp y_{[0]}
			& = e_{[-1]}  y_{[-2]} \sharp  (e_{[0]} \btl y_{[-1]} ) y_{[0]} \\
			& = e_{[-1]} y_{[-1]} \sharp y_{[0]} e_{[0]} \\
			& = (ye)_{[-1]} \sharp (ye)_{[0]}.
		\end{array}
		$$
			
	It remains to prove part (\ref{43}), that is 
		\begin{equation*}
			\lambda(\phi(f_{(1)} \btr x)) \cdot f_{(2)} = f \cdot \lambda(\phi(x)).
		\end{equation*}
	By equation (\ref{a}), this is 
		\begin{equation*}
			\lambda(\phi(x) \btl S(f_{(1)})) \cdot f_{(2)} = f \cdot \lambda(\phi(x)),
		\end{equation*}
	which is easily seen to follow from the right-left Yetter--Drinfeld condition (\ref{eq:YDrl}) for $\phi(x)$ and $Sf$ 
		\begin{equation*}
			(Sf)_{(2)} \cdot \lambda(\phi(x) \btl (Sf)_{(1)}) = \lambda(\phi(x)) \cdot Sf,			
		\end{equation*}
	that is 
		\begin{equation*}
			S(f_{(1)}) \cdot \lambda(\phi(x) \btl S(f_{(2)})) = \lambda(\phi(x)) \cdot S(f),
		\end{equation*} 
	by multiplying the latter on the right and on the left with certain components of $\Delta(f)$. With this, it is proven that $\Psi$ is a homomorphism. 
			
	Now we prove that it has a right inverse. By applying $S\otimes \phi^{-1}$ to formula (\ref{b}), we get 
		\begin{equation*} 
		S(\phi(x)_{[-1]}) \otimes \phi^{-1}(\phi(x)_{[0]}) = S^2(x_{[1]}) \otimes  x_{[0]}
		\end{equation*}
	and therefore
		\begin{equation*}
			f S(y_{[-1]}) \cdot \phi^{-1}( y_{[0]}) = fS^2({x_{[1]}}) \cdot x_{[0]},
		\end{equation*}
	where $x = \phi^{-1}(y)$. With the  homomorphism property of $\Psi$ proven and having that, by compatibility of coactions (\ref{b}), 
		\begin{equation*}  
			\Psi(x\sharp f) = S(x_{[1]}) \sharp \phi(x_{[0]}) \cdot f, 
		\end{equation*} 
	it is easy to check that $\Psi(f S^2(x_{[1]}) \cdot x_{[0]}) = f\sharp \phi(x)$. Therefore, its right inverse is 
		\begin{equation*} 
			f \sharp y \mapsto   f S^2(x_{[1]}) \cdot x_{[0]} =  f S(y_{[-1]}) \cdot \phi^{-1}( y_{[0]}),
		\end{equation*}
	where $x$ denotes the element of $L$ such that $\phi(x) = y$. To prove that it is also a left inverse, we first prove that $\Psi^{-1}$ is a homomorphism. From this it easily follows that
		\begin{align*}
			\Psi^{-1}(\Psi(x\sharp f) ) &= \Psi^{-1}(\Psi(x) \cdot \Psi( f) ) & \\
			& = \Psi^{-1}(\Psi(x)) \cdot \Psi^{-1}(\Psi( f)) & \\
			&= \Psi^{-1}(S(x_{[1]}) \cdot \phi(x_{[0]}) )\cdot \Psi^{-1}(f) & 
			\\
			&= S(x_{[1]}) \cdot  S(\phi(x_{[0]})_{[-1]}) \cdot \phi^{-1}(\phi(x_{[0]})_{[0]}) \cdot f & 
			\\
			&= S(x_{[2]}) S^2(x_{[1]}) \cdot x_{[0]} \cdot f & \text{by (\ref{b})} & \\
			&= x\sharp f. &
		\end{align*}
	Similarly as in Theorem \ref{BM}, after we prove the statement for special cases, for $y,c \in R$ and $f,g \in H$, 
	\begin{enumerate}[label={(\alph*)}, ref=\alph*]
		\item $\Psi^{-1}( f \cdot g) = \Psi^{-1}(f) \Psi^{-1}(g)$ \label{51}
		\item $\Psi^{-1}(f \sharp y)  = \Psi^{-1}(f) \Psi^{-1}(y)$ \label{52}
		\item $\Psi^{-1}(1\sharp y \cdot f\sharp 1) = \Psi^{-1}(y) \Psi^{-1}(f)$ \label{53}
		\item $\Psi^{-1}(y \cdot c) = \Psi^{-1}(y) \Psi^{-1}(c)$, \label{54}
	\end{enumerate}
	the homomorphism property follows directly. Parts (\ref{51}) and (\ref{52}) follow directly from the definition. Parts (\ref{53}) and (\ref{54}) are proven similarly as the corresponding equations in the proof of the antihomomorphism property of $\tau$ in Theorem \ref{BM}, as follows.
			
	We prove equation (\ref{53}). On the left-hand side we have
		\begin{align}
			\Psi^{-1}(y \cdot f) &=  \Psi^{-1}(f_{(1)} \sharp( y \btl f_{(2)})) \nonumber \\
			&= f_{(1)} S((y \btl f_{(2)})_{[-1]}) \cdot \phi^{-1}((y \btl f_{(2)})_{[0]}) \label{lhs} 
		\end{align} 
	and on the right-hand side
		\begin{align}
			\Psi^{-1}(y) \cdot \Psi^{-1}(f) &=  S(y_{[-1]}) \cdot \phi^{-1}(y_{[0]}) \cdot f. \label{rhs}
		\end{align} 
	Similarly as in Theorem \ref{BM}, we start from the right-left Yetter--Drinfeld property (\ref{eq:YDrl1}) and write consequently
		\begin{align*}
			f_{(1)} \otimes f_{(3)} (y\btl f_{(2)})_{[-1]} \otimes f_{(4)} \otimes (y\btl f_{(2)})_{[0]} \\ 
			= f_{(1)} \otimes y_{[-1]} f_{(2)} \otimes f_{(4)} \otimes (y_{[0]} \btl f_{(3)})
		\end{align*}
	and then apply $\id \otimes S \otimes \id \otimes \phi^{-1} $ and multiply tensor factors to obtain equality involving expression in  (\ref{lhs}),
		\begin{align}
			f_{(1)} S((y\btl f_{(2)})_{[-1]} ) \cdot \phi^{-1}((y\btl f_{(2)})_{[0]}) \nonumber \\ = S(y_{[-1]} )  f_{(2)} \cdot \phi^{-1}(y_{[0]} \btl f_{(1)}). \label{lhs2}
		\end{align}
	The equality of expressions in (\ref{lhs2}) and (\ref{rhs}), therefore (\ref{53}), follows from 
		\begin{align*}
			f_{(2)} \cdot \phi^{-1}(z \btl f_{(1)}) & =  (f_{(2)} \btr \phi^{-1}(z \btl f_{(1)}) ) \sharp f_{(3)} & \\
			&= \phi^{-1}((z \btl f_{(1)}) \btl S(f_{(2)})) ) \sharp f_{(3)} & \text{by (\ref{a})} \\
			& = \phi^{-1}(z) \sharp f &
		\end{align*}
	being true for any $z\in R$.	
		
	We now prove (\ref{54}). We have
		\begin{align*}
			\Psi^{-1}(y \cdot c) &= S((yc)_{[-1]}) \cdot \phi^{-1}((yc)_{[0]}) \\
			&= S(y_{[-1]})S(c_{[-1]}) \cdot \phi^{-1}(c_{[0]})\phi^{-1}(y_{[0]}) \\
			&= S(y_{[-1]}) \cdot \Psi^{-1}(c) \cdot \phi^{-1}(y_{[0]}), \\
			\Psi^{-1}(y) \cdot \Psi^{-1}(c) &= S(y_{[-1]}) \cdot \phi^{-1}(y_{[0]}) \cdot \Psi^{-1}(c).
		\end{align*}
	The equality follows from the fact that $$\Psi^{-1}(c) = S(c_{[-1]}) \cdot \phi^{-1}(c_{[0]}) = S^2(d_{[1]} )\cdot d_{[0]},$$ where $d = \phi^{-1}(c)$, commutes with every $x = x\sharp 1 \in L\sharp 1$, which is actually the auxiliary statement (\ref{aux}). It follows, as shown in the proof of Theorem \ref{BM}, from the braided commutativity of $(L,\btr,\rho)$. Therefore, it is proven that $\Psi$ is a homomorphism with inverse $\Psi^{-1}$.
			
	It is evident that $\Psi(L\sharp 1) = \operatorname{Im}\lambda$. We prove that $\Psi(\operatorname{Im} \rho) = 1\sharp R$. First, 
		\begin{align*}
			\Psi(\rho(x)) &= \Psi(x_{[0]} \sharp x_{[1]}) = \phi(x_{[0]})_{[-1]} \sharp \phi(x_{[0]})_{[0]} \cdot x_{[1]} \\
			& = S(x_{[1]}) \sharp \phi(x_{[0]}) \cdot x_{[2]} \\
			&= 1\sharp (\phi(x_{[0]}) \btl x_{[1]}) \in 1\sharp R
		\end{align*} 
	by compatibility of coactions. Therefore, $\Psi(\operatorname{Im}\rho) \subset 1\sharp R$. Next, 
		\begin{align*}
			\Psi^{-1}(1\sharp y) & = S(y_{[-1]}) \cdot \phi^{-1}(y_{[0]}) = S^2(x_{[1]}) \cdot x_{[0]},
		\end{align*}  
	where $x = \phi^{-1}(y)$. This is in $\operatorname{Im}\rho$ because we have 
		\begin{align}
			\rho(S^2(x_{[1]}) \btr x_{[0]}) & = \rho(S^2(x_{[3]}) \btr x_{[0]}) \cdot  S^2(x_{[2]}) \cdot S(x_{[1]}) \nonumber \\ 
			& = S^2(x_{[2]}) \cdot \rho(x_{[0]}) \cdot S(x_{[1]}) \nonumber\\
			& = S^2(x_{[3]}) \cdot x_{[0]} \sharp x_{[1]} \cdot S(x_{[2]}) \nonumber \\
			& = S^2(x_{[1]}) \cdot x_{[0]}, \label{Imrho}
		\end{align}
	by using the left-right Yetter--Drinfeld property (\ref{eq:YD}) for $x_{[0]}$ and $S^2(x_{[2]})$ in the second row.

	We prove (\ref{33}). Now it is easy to check that $ \Psi \circ \tau  = \Phi$ and $\tau' \circ \Psi = \Phi$.
		\begin{align*}
			\Psi(\tau(x\sharp f)) & = \Psi(SfS^{2}( x_{[1]}) \cdot x_{[0]}) & \\
			& = \Psi(SfS^{2}(x_{[1]})) \cdot \Psi(x_{[0]}) & \text{by (\ref{32})} \\
			& = SfS^{2}(x_{[2]}) \cdot S(x_{[1]})\sharp \phi(x_{[0]}) &  \\
			& = Sf \sharp \phi(x)& \\
			& = \Phi(x\sharp f), & \\
			\Phi(\Psi^{-1}(f\sharp y)) &= \Phi(fS(y_{[-1]}) \cdot \phi^{-1}(y_{[0]})) & \\
			&= \Phi(\phi^{-1}(y_{[0]}) \cdot \Phi(S(y_{[-1]})) \cdot  \Phi(f)  & \text{by (\ref{31})}\\
			&= y_{[0]} \cdot S^2(y_{[-1]}) Sf \\
			&= \tau'(f\sharp y).
		\end{align*}
	Since $\Phi$ is an antihomomorphism and $\Psi$ an isomorphism, it follows that $\tau = \Psi^{-1} \circ \Phi$ and $\tau' = \Phi \circ \Psi^{-1}$ are antihomomorphisms. 
		
\end{proof}

Now we can transfer the structure of a right bialgebroid from $H\sharp R$ to $L\sharp H$, and the structure of a left bialgebroid from $L\sharp H$ to $H\sharp R$, by using the algebra iso\-mor\-phism~$\Psi$, and prove that these comprise the structure of a symmetric Hopf algebroid on $L\sharp H \cong H\sharp R$, with antipode $\tau$.
	
\begin{theorem} \label{HAHA}
	Let $H$ be a Hopf algebra with an antipode $S$. Let $\phi \colon L \to R$ be an antiisomorphism of  algebras. Let $(L, \btr, \rho)$ and $(R,\btl, \lambda)$ a left-right and a right-left braided commutative Yetter--Drinfeld $H$-module algebras compatible via $\phi$, that is, for all $x\in L$ and $f \in H$
		\begin{align*}
			\phi(f \btr x) &= \phi(x) \btl Sf, \\
			\phi(x)_{[-1]} \otimes \phi(x)_{[0]} &= S(x_{[1]}) \otimes \phi( x_{[0]})  .
		\end{align*}
	\begin{enumerate}
		\item Then $L\sharp H$ with the following structure maps is a symmetric Hopf algebroid. \label{61}
			\begin{align*} 
				\alpha_L \colon & L \to L\sharp H &  && \alpha_L(x) &= x \sharp 1_H \\
				\beta_L  \colon & L \to L \sharp H & && \beta_L(x) &= \rho(x) =  x_{[0]} \sharp  x_{[1]}\\
				\Delta_L  \colon & L \sharp H \to L\sharp H \otimes_L L \sharp H & && \Delta_L(x \sharp f) &= x \sharp f_{(1)} \otimes_L 1_L \sharp f_{(2)} \\
				\epsilon_L \colon & L  \sharp H \to L & && \epsilon_L(x\sharp f ) &= \epsilon(f) x \\
				\\
				\alpha_R  \colon& R \to L\sharp H &&& \alpha_R( y ) &= 
				S(y_{[-1]}) \cdot \phi^{-1}(y_{[0]})\\
				\beta_R  \colon & R \to L\sharp H &&& \beta_R( y ) &= \phi^{-1}(y) \sharp 1_H \\ 
				\Delta_R  \colon & L \sharp H \to L\sharp H \otimes_{R} L \sharp H &&& \Delta_R( x \sharp f ) &= x \sharp f_{(1)} \otimes_R 1_L \sharp f_{(2)}\\
				\epsilon_R \colon& L \sharp H \to R &&& \epsilon_R(x \sharp f) &= \phi(x) \btl f \\
				\\ 
				\tau  \colon & L \sharp H \to L \sharp H &&& \tau(x \sharp f) &= Sf S^2(x_{[1]}) \cdot x_{[0]}.
			\end{align*}
			
		\item Then $H\sharp R$ with the following structure maps is a symmetric Hopf algebroid. \label{62}
			\begin{align*} 
				\alpha_L' \colon& L \to H\sharp R &&& \alpha_L'(x) &= \lambda(\phi(x)) = S(x_{[1]}) \sharp  \phi(x_{[0]}) \\
				\beta_L'  \colon& L \to H\sharp R &&& \beta_L'(x) &= 
				\phi(x_{[0]})\btl x_{[1]}\\ 
				\Delta_L'  \colon & H \sharp R \to H\sharp R \otimes_{L} H \sharp R &&& \Delta_L'(f \sharp y) &= f_{(1)}\sharp 1_R \otimes_L f_{(2)}\sharp y \\
				\epsilon_L' \colon& H \sharp R \to L &&& \epsilon_L'(f\sharp y) &= 
				f S(y_{[-1]}) \btr \phi^{-1}(y_{[0]}) 			\\
				\\
				\alpha_R'  \colon& R \to H\sharp R &&&  \alpha_R'(y) &= 1_H \sharp y \\
				\beta_R'  \colon& R \to H\sharp R &&& \beta_R'(y) &= \lambda(y) =  y_{[-1]} \sharp  y_{[0]}\\
				\Delta_R'  \colon & H \sharp R \to H\sharp R \otimes_R H \sharp R  &&& \Delta_R'( f \sharp y) &= f_{(1)}\sharp 1_R \otimes_R f_{(2)}\sharp y \\
				\epsilon_R' \colon& H \sharp R \to R &&& \epsilon_R'(f\sharp y) &= \epsilon(f) y \\
				\\
				\tau'  \colon & H \sharp R \to H \sharp R &&& \tau'(f \sharp y ) & = y_{[0]} \cdot S^2(y_{[-1]})Sf.
			\end{align*}
	\end{enumerate}
	These are two smash product algebra presentations of the same symmetric Hopf algebroid, connected by isomorphism $\Psi \colon L\sharp H \to H\sharp R$ defined in Theorem \ref{Phi}. Algebras $H$, $L$ and $R$ are canonically embedded as subalgebras in algebra $L\sharp H \cong H\sharp R$. Furthermore, subalgebra $L$ is the image of the canonical embedding $\alpha_L$ and the image of~$\beta_R$, and subalgebra $R$ is the image of the canonical embedding $\alpha_R$ and the image of~$\beta_L$. The antiisomorphism $\phi$ satisfies $\phi = \epsilon_R \circ \alpha_L$ and its inverse is $\phi^{-1} = \epsilon_L \circ \beta_R$.
\end{theorem}
	
\begin{proof}
	 We use notation for structure maps from Proposition \ref{bial} and isomorphism $\Psi \colon L\sharp H \to H \sharp R$ from Theorem \ref{Phi}, with formulas
		\begin{align*}
		\Psi ( x \sharp f ) & = S(x_{[1]}) \sharp \phi(x_{[0]}) \cdot f, \\
		\Psi^{-1}(f \sharp y) & =  f S(y_{[-1]}) \cdot \phi^{-1}( y_{[0]}) = f S^2(x_{[1]}) \cdot x_{[0]},
		\end{align*}
	where $x = \phi^{-1}(y) $, $y\in R$ and $f \in H$.
	It is easy to see that source and target maps agree, by checking $$ \alpha_R = \Psi^{-1} \circ \alpha_R',\quad \beta_R = \Psi^{-1} \circ \beta_R',\quad \Psi \circ \alpha_L = \alpha_L',\quad \Psi \circ \beta_L = \beta_L',$$ and that counits agree, by checking $$\epsilon_R= \epsilon_R'  \circ \Psi,\quad \epsilon_L \circ \Psi^{-1}= \epsilon_L'.$$ We note that this puts the $L$-bimodule structure on $H\sharp R$ and that $\Psi \otimes \Psi$ is well defined as a map from $L\sharp H \otimes_L L\sharp H $ to $H\sharp R \otimes_L H\sharp R$, because $\Psi$ is a homomorphism and source and target maps agree. Similar is true for $L\sharp H$, which now has also a structure of an $R$-bimodule. 
		
	Now we check that coproducts agree, that is $$(\Psi \otimes \Psi )\circ \Delta_L = \Delta_L' \circ \Psi,\quad (\Psi \otimes \Psi) \circ \Delta_R = \Delta_R' \circ \Psi.$$ We use that $H\sharp R$ is an $L$-bimodule with $x\in L$ acting on the left by multiplication on the left with $\alpha_L'(x)$, and $x$ acting on the right by multiplication on the left with $\beta_L'(x)$, to prove the first equality,
		\begin{align*}
		(\Psi\otimes\Psi)(\Delta_L(x\sharp f)) &= \Psi (x\sharp f_{(1)})\otimes_L \Psi(f_{(2)}) \\
		&= S(x_{[1]}) \sharp \phi(x_{[0]}) \cdot f_{(1)} \otimes_L f_{(2)}, \\
		 \Delta_L' ( \Psi(x\sharp f)) &= \Delta_L' ( S(x_{[1]}) \sharp \phi(x_{[0]}) \cdot f_{(1)}) \\
		&= \Delta_L' ( S(x_{[1]}) f_{(1)} \sharp ( \phi(x_{[0]}) \btl f_{(2)})) \\
		&=  S(x_{[2]}) f_{(1)} \otimes_L S(x_{[1]}) f_{(2)} \sharp ( \phi(x_{[0]}) \btl f_{(3)}) \\
		&=  S(x_{[2]}) f_{(1)} \otimes_L S(x_{[1]}) \sharp  \phi(x_{[0]}) \cdot f_{(2)} \\
		&=  S(x_{[1]}) f_{(1)} \otimes_L \alpha_L'(x_{[0]}) \cdot f_{(2)}\\
		&= \beta_L'(x_{[0]})  \cdot S(x_{[1]}) f_{(1)}  \otimes_L  f_{(2)} \\
		&= (\phi(x_{[0]}) \btl x_{[1]}) \cdot S(x_{[2]}) f_{(1)}  \otimes_L  f_{(2)} \\
		&= S(x_{[1]}) \sharp \phi(x_{[0]}) \cdot f_{(1)} \otimes_L f_{(2)} .
		\end{align*}
	We use the right-left Yetter--Drinfeld condition to prove the second equality,
		\begin{align*}
		(\Psi\otimes\Psi)(\Delta_R(x\sharp f)) &= \Psi (x\sharp f_{(1)})\otimes_R \Psi(f_{(2)}) \\
		&= S(x_{[1]}) \sharp \phi(x_{[0]}) \cdot f_{(1)} \otimes_R f_{(2)} \\
		&= S(x_{[1]}) f_{(1)} \sharp (\phi(x_{[0]}) \btl f_{(2)}) \otimes_R f_{(3)} \\
		&= S(x_{[1]}) f_{(1)} \cdot \alpha'_R(\phi(x_{[0]}) \btl f_{(2)}) \otimes_R f_{(3)} \\
		&= S(x_{[1]}) f_{(1)} \otimes_R f_{(3)} \cdot \beta'_R(\phi(x_{[0]}) \btl f_{(2)})  \\
		&= S(x_{[1]}) f_{(1)} \otimes_R f_{(3)} \cdot \lambda(\phi(x_{[0]}) \btl f_{(2)})  \\
		&= S(x_{[1]}) f_{(1)} \otimes_R   \lambda(\phi(x_{[0]})) \cdot f_{(2)},  \\
		\Delta_R' ( \Psi(x\sharp f)) &=  \Delta_R' ( S(x_{[1]}) f_{(1)} \sharp ( \phi(x_{[0]}) \btl f_{(2)})) \\
		&=  S(x_{[2]}) f_{(1)} \otimes_R S(x_{[1]}) f_{(2)} \sharp ( \phi(x_{[0]}) \btl f_{(3)}) \\
		&=  S(x_{[2]}) f_{(1)} \otimes_R S(x_{[1]}) \sharp  \phi(x_{[0]}) \cdot f_{(2)} .
		\end{align*}
		
	Since, by Proposition \ref{bial}, $H\sharp R$ is a right bialgebroid and $L\sharp H$ a left bialgebroid, and  all structure maps here agree with structure maps of these scalar extension bialgebroids, it follows that $L\sharp H$ is a right bialgebroid and $H\sharp R$ is a left bialgebroid. By Theorem \ref{Phi}, the maps $\tau$ and $\tau'$ are antihomomorphisms and they satisfy equality $\Psi \circ \tau = \tau' \circ \Psi$.
		
	 Formulas for compatibilities of Hopf algebroid structure maps
		\begin{equation*} 
			\begin{aligned} 
			\alpha_L\circ\epsilon_L\circ\beta_R = \beta_R,
			&\quad
			\beta_L\circ\epsilon_L\circ\alpha_R = \alpha_R, 
			\\
			\alpha_R\circ\epsilon_R\circ\beta_L = \beta_L,
			&\quad
			\beta_R\circ\epsilon_R\circ\alpha_L = \alpha_L,  \\
			(\Delta_R\otimes_L\id)\circ\Delta_L &=
			(\id\otimes_R\Delta_L)\circ\Delta_R,   \\
			(\Delta_L\otimes_R\id)\circ\Delta_R &=
			(\id\otimes_L \Delta_R)\circ\Delta_L, \\
			\tau\circ\beta_L =\alpha_L, &\quad
			\tau\circ\beta_R = \alpha_R, \\
			\mu_{\otimes'_L}\circ(\tau\otimes\id)&\circ\Delta_L = \alpha_R\circ\epsilon_R, 
			\\
			\mu_{\otimes'_R}\circ(\id\otimes \tau)&\circ\Delta_R = \alpha_L\circ\epsilon_L ,
			\end{aligned}
		\end{equation*}
	are all very easy to prove when we chose to prove each in the smash product algebra that has simpler formulas for the corresponding structure maps, with the exception of the second formula. In more detail, the first formula and the fourth formula are $\alpha_L \circ \phi^{-1} = \beta_R$ and $ \beta_R \circ \phi = \alpha_L$, respectively. The third formula $\alpha'_R\circ\epsilon'_R\circ\beta'_L = \beta'_L$ is $\phi(x_{[0]}) \btl x_{[1]} = \phi(x_{[0]}) \btl x_{[1]}$. The second formula $\beta_L\circ\epsilon_L\circ\alpha_R = \alpha_R$ is 
		\begin{equation*}
			(S^2(x_{[1]}) \btr x_{[0]}  )_{[0]}\sharp  (S^2(x_{[1]}) \btr x_{[0]}  )_{[1]} = S^2(x_{[1]}) \cdot x_{[0]},
		\end{equation*}
	which is proven in (\ref{Imrho}) inside the proof of Theorem \ref{Phi}, by using the left-right Yetter--Drinfeld condition for $(L,\btr,\rho)$. Next, formula $\tau \circ \beta_L = \alpha_L$ is true equation $S(x_{[2]}) S^2(x_{[1]}) \cdot x_{[0]} = x$, and formula $\tau \circ \beta_R' = \alpha_R'$ holds similarly. 
	Last two formulas are respectively $S(f_{(1)}) \cdot f_{(2)} \sharp y = y$ in $H\sharp R$ and $x \sharp f_{(1)} \cdot S(f_{(2)}) = x$ in $L\sharp H$.
		
\end{proof}


\subsection{Reasons for the compatibility. Another formulation} \label{sectionB}
	
Given any symmetric Hopf algebroid structure, there are two canonical antiisomorphisms of $L$ and $R$ that arise \cite{bohmHbk}. These are $\phi\colon L \to R$, $\phi = \epsilon_R \circ \alpha_L$, with inverse $\phi^{-1} = \epsilon_L \circ \beta_R$, and symmetrically $\theta\colon R \to L$, $\theta = \epsilon_L \circ \alpha_R$, with inverse $\theta^{-1} = \epsilon_R \circ \beta_L$. In the previous subsection, we started with any given antiisomorphism $\phi \colon L  \to R$ and two braided commutative Yetter--Drinfeld module algebras compatible via $\phi$ and obtained a symmetric Hopf algebroid in which $\phi$ is equal to $\epsilon_R \circ \alpha_L$. The analogous construction can be done with antiisomorphism $\theta$, which is shown below. First, we show how the compatibility of Yetter--Drinfeld module algebras arises naturally when considering an algebra isomorphism of one left and one right scalar extension bialgebroid. In this way, we then obtain the definition of the compatibility via $\theta$ and the statement of the corresponding analogous theorem.
	
\subsubsection{Right bialgebroid with additional compatible right coaction} 
	
Starting with a right-left braided commutative Yetter--Drinfeld $H$-module algebra $(R,\btl,\lambda)$, that is, a right bialgebroid $(H\sharp R, \alpha'_R, \beta'_R, \Delta'_R, \epsilon'_R)$, the image of the injective antihomomorphism $\beta'_R$ can be viewed as a subalgebra $L \subset H\sharp R$, which we would like to use to present the algebra $H\sharp R$ as another smash product algebra $L\sharp H$. If this another presentation exists, it is unique, because the corresponding left Hopf action is uniquely determined by the multiplication on the algebra $H\sharp R$, somewhat unprecisely written as
	\begin{equation}\label{smasha}
	f \btr x = f_{(1)} \cdot x \cdot S(f_{(2)}), \quad \text{ for } f\in H \text{ and } x\in L.
	\end{equation} 
Let now $\alpha'_L \colon L \to H\sharp R$ be the injective homomorphism with image $\operatorname{Im} \beta'_R$ formalizing this inclusion $L\subset H\sharp R$. Denote by $\phi$ the antiisomorphism such that $\beta'_R \circ \phi = \alpha'_L$, that is 
	\begin{equation}\label{alpha}
	\alpha'_L(x) = \beta'_R(\phi(x)) = \phi(x)_{[-1]} \otimes \phi(x)_{[0]}.
	\end{equation} 
We apply $\epsilon'_R$ to this and get $\epsilon'_R \circ \alpha'_L  = \phi$, a formula for $\phi$. Now (\ref{smasha}) can be written more precisely as 
	\begin{equation}\label{pre}
	\alpha'_L(f \btr x) = f_{(1)} \cdot \phi(x)_{[-1]} \sharp \phi(x)_{[0]} \cdot S(f_{(2)}).
	\end{equation}
When we apply $\epsilon'_R$ to it, we get a formula that the left action should satisfy:
	\begin{equation}
	\phi(f \btr x)= \phi(x) \btl Sf.
	\end{equation} 
We now use this formula to define this left action. It can be easily proven that it is a left Hopf action, thus with the initial right bialgebroid $H\sharp R$ we naturally get a smash product algebra $L\sharp H$. The map 
	\begin{equation} \label{psialpha}
	\Psi\colon L\sharp H \to H \sharp R, \quad \Psi(x\sharp f) = \alpha'_L(x) \cdot f = \beta'_R(\phi(x))\cdot f
	\end{equation} 
is automatically a surjective homomorphism with $\Psi(L\sharp 1) = \operatorname{Im}\beta_R'$, which we now prove. First, it is a homomorphism if 
	$
	\Psi(f \cdot x) = \Psi(f) \cdot \Psi(x),
	$
since other cases can be easily checked to be true. We have
	\begin{align*}
	\Psi(f \cdot x) & = \Psi((f_{(1)} \btr x) f_{(2)}) \\
	& = \phi(f_{(1)} \btr x)_{[-1]} \sharp \phi(f_{(1)} \btr x)_{[0]} \cdot f_{(2)} \\
	& = (\phi(x) \btl S(f_{(1)}))_{[-1]} \sharp (\phi(x) \btl S(f_{(1)}))_{[0]} \cdot f_{(2)}, \\
	\Psi(f) \cdot \Psi(x) &= f  \phi(x)_{[-1]} \sharp \phi(x)_{[0]}.
	\end{align*}
The equality
	\begin{equation}
	(\phi(x) \btl Sf)_{[-1]} \sharp (\phi(x) \btl Sf)_{[0]}  = f_{(1)}  \phi(x)_{[-1]} \sharp \phi(x)_{[0]} \cdot S(f_{(2)})
	\end{equation}
is  equivalent to the right-left Yetter--Drinfeld condition for $Sf$ and $\phi(x)$,
	\begin{equation}
	(Sf)_{(2)}(\phi(x) \btl (Sf)_{(1)})_{[-1]} \sharp (\phi(x) \btl (Sf)_{(1)})_{[0]}  =  \phi(x)_{[-1]} \smash \phi(x)_{[0]} \cdot Sf,
	\end{equation}
hence it holds and $\Psi$ is a homomorphism. Second, from this, it follows that the preimages of elements $ f\sharp 1 \in H \sharp R$ and $1\sharp y \in H\sharp R$ exist: they are easily checked to be $1\sharp f$ and $S(y_{[-1]}) \cdot \phi^{-1}(y_{[0]}) \in L\sharp H$, respectively. Together with the homomorphism property, this implies that $\Psi$ is surjective. 
	
Therefore, this new smash product algebra is a different presentation of the same algebra if and only if $\Psi$ is  injective. 

The antipode $\tau' \colon H\sharp R \to H \sharp R$, if it exists, is uniquely determined because it is an antihomomorphism that is on subalgebra $H$ equal to $S$ and on subalgebra $L = \Psi(L\sharp 1) \subset  H\sharp R$ equal to the left inverse of $\beta'_R$, hence
	\begin{equation}
	\tau'(f \sharp y )  = y_{[0]} \cdot S^2(y_{[-1]})Sf. \label{ttau}
	\end{equation} 
On the other hand, the map defined by (\ref{ttau}) is automatically an antihomomorphism, which is for $\tau'$ proven analogously to proof for $\tau$ in Theorem \ref{BM}. Therefore, this antipode always exists and is unique. 
	
Now, assume that we have additionally a right comodule structure $(L,\rho)$ such that $\rho \colon L \to L\sharp H$ is an antihomomorphism with $\Psi(\operatorname{Im}\rho) = 1\sharp R$, and set $\beta'_L = \Psi \circ \rho$. For  $\alpha'_L = \tau' \circ \beta'_L$ to be true, this coaction $\rho \colon x \mapsto x_{[0]} \otimes x_{[1]}$ has to be compatible via $\phi$ with the coaction $\lambda$,
	\begin{equation}
	\phi(x)_{[-1]} \otimes \phi(x)_{[0]} = S(x_{[1]}) \otimes \phi( x_{[0]}).
	\end{equation}
This follows from $\alpha'_L = \tau' \circ \beta'_L = \tau' \circ \Psi \circ \rho$, 	since the antihomomorphism $\Phi := \tau' \circ \Psi$ equals by a simple calculation 
	\begin{equation}
	\Phi(x\sharp f) = Sf \sharp \phi(x).
	\end{equation}
On the other hand, it is easy to see that the compatibility of coactions implies $\Psi(\operatorname{Im}\rho) = 1\sharp R$ and $\alpha_L' = \tau' \circ \beta_L'$.
	
In the proof of Theorem \ref{Phi}, we see that injectivity of $\Psi$ follows from the antihomomorphism property of $\phi$, $\lambda$ and $\rho$, the compatibility of actions and coactions via~$\phi$ and the right-left Yetter--Drinfeld condition on $(R,\btl,\lambda)$. In Theorem \ref{obrat} below, it is proven that isomorphism property of $\Psi$ and $\Psi(\operatorname{Im}\rho) = 1\sharp R$ imply the left-right Yetter--Drinfeld property for $(L,\btr, \rho)$, and therefore a symmetric Hopf algebroid structure by Theorem~\ref{HAHA}.
  
In short, starting with a right scalar extension bialgebroid $H\sharp R$, the only ad\-di\-ti\-o\-nal data we need is a comodule structure $(L,\rho)$ such that $\rho \colon L\to L\sharp H$ is an antihomomorphism compatible via $\phi$ with coaction $\lambda$. Here the action that gives rise to the codomain smash product algebra $L\sharp H$ is defined by $f\btr x = \phi^{-1}(\phi(x) \btl Sf)$, where antiisomorphism $\phi$ is defined as in the beginning. 
	
Thus we have shown how antiisomorphism $\phi$, isomorphism $\Psi$ and compatibility of actions and coactions via $\phi$ arise naturally when considering possible compatible left bialgebroid structure on the right  bialgebroid $H \sharp R$.

\subsubsection{Left bialgebroid with additional compatible left coaction}

We briefly walk through the analogous discussion for a left bialgebroid, to obtain the statement of the theorem with aforementioned $\theta$ analogous to the statement of Theorem \ref{Phi}. 
	
Starting with a left-right braided commutative Yetter--Drinfeld $H$-module algebra $(L,\btr,\rho)$, that is, a left bialgebroid $(L\sharp H,\alpha_L,\beta_L,\Delta_L,\epsilon_L)$, the image of the injective antihomomorphism $\beta_L$ can be viewed as a subalgebra $R \subset L\sharp H$. Let  $\alpha_R \colon R \to L\sharp H$ be the injective homomorphism with the image $\operatorname{Im}\beta_L$ formalizing this inclusion. This, similarly as before, can possibly give rise to a different presentation $H\sharp R$. The corresponding right action is uniquely determined by the multiplication on the algebra~$L\sharp H$, for $y \in R$ and $f \in H$,
$$y \btl f = S(f_{(1)}) \cdot y \cdot f_{(2)}, $$ which can be more precisely written as
	\begin{equation}\label{amu}
	\theta(y \btl f) = Sf \btr \theta(y),
	\end{equation} 
with antiisomorphism $\theta \colon R \to L$ defined by $\beta_L \circ \theta = \alpha_R$, hence $\theta= \epsilon_L \circ \alpha_R$. The surjective homomorphism between these two algebras is $\bar \Psi \colon H\sharp R \to L \sharp H$,
	\begin{equation}
	\bar \Psi(f\sharp y ) = f \cdot \alpha_R(y) = f \cdot \theta(y)_{[0]} \sharp \theta(y)_{[1]},
	\end{equation}
hence they are isomorphic if and only if $\bar \Psi$ is injective. 
The antipode $\tau$ is uniquely determined since it is an antihomomorphism that is on subalgebra $H$ equal to $S$ and on subalgebra $R = \bar\Psi(1\sharp R) \subset L\sharp H$ equal to the left inverse of $\beta_L$,
	\begin{equation}
	\tau(x\sharp f)  = Sf S^2(x_{[1]}) \cdot x_{[0]}.
	\end{equation} 
Assume now we have additionally a left comodule structure $(R,\lambda)$ such that map $\lambda \colon R \to H\sharp R$ is an antihomomorphism with $\bar\Psi(\operatorname{Im}\lambda) = L\sharp 1$, and set $\beta_R := \bar \Psi \circ \lambda$. For $\alpha_R = \tau \circ \beta_R$ to hold, coaction $\lambda$ has to be compatible with $\rho$ as in
	\begin{equation}\label{bmu}
	\theta(y)_{[0]} \otimes \theta(y)_{[1]} =  \theta( y_{[0]}) \otimes S(y_{[-1]}).
	\end{equation}
This follows from $\alpha_R = \tau \circ \beta_R = \tau \circ \bar \Psi \circ \lambda$ by easy calculation $$\tau(\bar\Psi(f\sharp y)) = \theta(y) \sharp Sf.$$

\subsubsection{Analogous theorem with another compatibility}
	
Previous theorem and this discussion together point to the analogous theorem with antiisomorphism $\theta$ and the compatibility defined by formulas (\ref{amu}) and (\ref{bmu}), which we give here without proof.
	
\begin{definition} \label{comptheta} 
	Let $H$ be a Hopf algebra with an antipode $S$. Let $\theta \colon R \to L$ be an antiisomorphism of algebras. We say that a left-right braided commutative Yetter--Drinfeld $H$-module algebra $(L,\btr,\rho)$ and a right-left braided commutative Yetter--Drin\-feld $H$-module algebra $(R,\btl,\lambda)$ are \emph{compatible via} $\theta$ if their structure maps satisfy, for all $y \in R$ and $f \in H$, 
		\begin{align*}
		\theta( y \btl f) &=Sf \btr \theta(y) ,\\ 
		\theta(y)_{[0]} \otimes \theta(y)_{[1]} &= \theta (y_{[0]}) \otimes S( y_{[-1]}).
		\end{align*} 
\end{definition}

\begin{theorem} \label{barPhi}
	Let $H$ be a Hopf algebra with an antipode $S$. Let $\theta \colon R \to L$ be an antiisomorphism of  algebras. Let $(L, \btr, \rho)$ and $(R,\btl, \lambda)$ be a left-right and a right-left braided commutative Yetter--Drinfeld $H$-module algebras compatible via $\theta$, that is, for all $y \in R$ and $f \in H$, 
		\begin{align}
		\theta(y \btl f) &= Sf \btr \theta(y), \label{ax}
		\\
		\theta(y)_{[0]} \otimes \theta(y)_{[1]} &= \theta(y_{[0]}) \otimes S( y_{[-1]}). \label{bx} 
		\end{align} 
	Then the map $\bar\Psi \colon H\sharp R \to L \sharp H$ defined by $$\bar\Psi(f\sharp y)=  f \cdot \theta(y)_{[0]}  \sharp \theta(y)_{[1]}$$ is an algebra isomorphism with inverse $\bar\Psi^{-1} \colon L\sharp H \to H \sharp R$ given by  $$\bar\Psi^{-1} (x\sharp f) = \theta^{-1}(x_{[0]})\cdot  S(x_{[1]}) f.$$ Also, $\bar\Psi(1\sharp R) = \operatorname{Im} \rho$ and $\bar\Psi^{-1}(L\sharp 1) = \operatorname{Im}\lambda$. The map $\bar\Phi \colon H\sharp R \to L \sharp H$ defined by  $$\bar\Phi(f\sharp y) = \theta(y) \sharp Sf$$ is an antihomomorphism. The maps $\tau$ and $\tau'$ defined by formulas below satisfy $\tau \circ \bar\Psi = \bar\Psi \circ \tau' =  \bar\Phi$ and are therefore antihomomorphisms. 
		
	Furthermore, the algebra $H\sharp R \cong L\sharp H$ carries a structure of a symmetric Hopf algebroid, with left bialgebroid structure maps as in left scalar extension bialgebroid when written in terms of the smash product algebra $L\sharp H$, right bialgebroid structure maps as in right scalar extension bialgebroid when written in terms of the smash product algebra $H\sharp R$, and antipode $\tau$ as below. 		
	The formulas for the structure maps written in terms of the other smash product algebra are also given below.
		\begin{align*}
		& L\sharp H&        &             H \sharp R  \\
		& \text{left bialgebroid over } L  & &  \text{left bialgebroid over }  L \\ 
		& \alpha_L(x)= x \sharp 1_H &              &                 \alpha'_L(x) =   \theta^{-1}(x_{[0]}) \cdot S(x_{[1]}) \\
		&\beta_L(x) = \rho(x) = x_{[0]} \sharp x_{[1]}          & & \beta'_L(x) =     1_H \sharp \theta^{-1}(x)   \\
		& \Delta_L(x \sharp f) = x \sharp f_{(1)} \otimes_L 1_L \sharp f_{(2)}      &  & \Delta'_L (f \sharp y) =  f_{(1)} \sharp 1_R \otimes_L f_{(2)} \sharp y  \\
		& \epsilon_L(x \sharp f) =\epsilon(f) x &   &   \epsilon'_L(f\sharp y)  =  f \btr \theta( y )  
		\\ \\
		& \text{right bialgebroid over } R & &  \text{right bialgebroid over }  R  \\   
		& \alpha_R(y)= \rho(\theta(y))=\theta(y_{[0]}) \sharp S(y_{[-1]}) & &                              \alpha'_R(y) =  1_H \sharp y   \\
		& \beta_R(y) =  y_{[-1]} \btr \theta(y_{[0]})       & & \beta'_R(y) = \lambda(y)= y_{[-1]} \sharp y_{[0]}  \\
		&  \Delta_R(x \sharp f) = x \sharp f_{(1)} \otimes_R 1_L \sharp f_{(2)}      &   & \Delta'_R(f \sharp y) = f_{(1)} \sharp 1_R \otimes_R f_{(2)} \sharp y   \\
		& \epsilon_R(x \sharp f) =\theta^{-1}(x_{[0]}) \btl S(x_{[1]})f  & &   \epsilon'_R(f \sharp y) = \epsilon(f)y   \\
		\\
		& \text{antipode}  &  & \text{antipode}  \\ 
		& \tau (x \sharp f) =  Sf S^2(x_{[1]}) \cdot x_{[0]} & &  \tau' (f \sharp y) = y_{[0]} \cdot S^2(y_{[-1]})Sf .
		\end{align*}
	Equations $\theta = \epsilon_L \circ \alpha_R = \epsilon'_L \circ \alpha'_R$ and $\theta^{-1} = \epsilon_R \circ \beta_L = \epsilon'_R \circ \beta'_L$ hold.	
\end{theorem}

\begin{proof}
	The proof is analogous to the proofs of Theorem \ref{Phi} and Theorem \ref{HAHA} together.
\end{proof}

\subsubsection{Equality of the two symmetric Hopf algebroids}

Let $(L,\btr,\rho)$ and $(R,\btl, \lambda)$ be braided commutative Yetter--Drinfeld module algebras compatible via $\theta$ as in the assumptions of the previous Theorem \ref{barPhi}. Then, by this theorem, a symmetric Hopf algebroid structure is defined on $L\sharp H \cong H\sharp R$ via isomorphism $\bar\Psi$. 
Now, for a map $\phi:= \epsilon_R \circ \alpha_L$ in this Hopf algebroid, one would expect that the Yetter--Drinfeld mo\-du\-le algebras are compatible also via $\phi$ and that the corresponding symmetric Hopf algebroid from Theorem \ref{HAHA}, defined by using isomorphism $\Psi$, is equal to the Hopf algebroid in the previous theorem, with $\bar\Psi = \Psi^{-1}$. This is true. It is proven in the next theorem.
	
\begin{theorem} \label{compa}
	Let $H$ be a Hopf algebra with an antipode $S$. 
	Let $(L, \btr, \rho)$ and $(R,\btl, \lambda)$ be a left-right and a right-left braided commutative Yetter--Drinfeld $H$-module algebras.
	\begin{enumerate} 
		
		\item Let $\phi \colon L \to R$ be an antiisomorphism of  algebras. If actions and coactions are compatible via $\phi$, that is, for all $f \in H$ and $x\in L$, \label{71}
			\begin{align*}
			\phi(f \btr x) &= \phi(x) \btl Sf ,
			\\
			\phi(x)_{[-1]} \otimes \phi(x)_{[0]} &= S(x_{[1]}) \otimes \phi( x_{[0]}), 
			\end{align*} 
		then actions and coactions are also compatible via $\theta$, 		
		where $\theta \colon R \to L$ is defined by  $\theta(y) = S(y_{[-1]}) \btr \phi^{-1}(y_{[0]})$.
	
		\item Let $\theta \colon R \to L$ be an antiisomorphism of  algebras. If actions and coactions are compatible via $\theta$, that is, for all $f \in H$ and $y\in R$, \label{72}
			\begin{align}
				\theta(y \btl f) &= Sf \btr \theta(y) , \label{aaa}
				\\
				\theta(y)_{[0]} \otimes \theta(y)_{[1]} &= \theta(y_{[0]}) \otimes S( y_{[-1]}), \label{bbb}
			\end{align} 
		then actions and coactions are also compatible via $\phi$,  where $\phi \colon L \to R$ is defined by  $\phi(x) = \theta^{-1}(x_{[0]}) \btl S(x_{[1]})$.
	\end{enumerate}
	When any of these two equivalent compatibilities of Yetter--Drinfeld module algebras holds, the resulting symmetric Hopf algebroids from Theorem \ref{HAHA} and Theorem \ref{barPhi} are the same, with $\theta = \epsilon_L \circ \alpha_R$ and $\phi = \epsilon_R \circ \alpha_L$.
\end{theorem}

\begin{proof} 
	We prove (\ref{72}). By Theorem \ref{barPhi}, we have a symmetric Hopf algebroid structure on $L\sharp H \cong H\sharp R$ via compatibility $\theta$. From formulas for structure maps, we get $\phi(x) = \theta^{-1}(x_{[0]}) \btl S(x_{[1]}) = \epsilon_R' (\alpha_L'(x))$ for all $x\in L$. 
	Therefore, indeed $\phi = \epsilon_R \circ \alpha_L = \epsilon'_R \circ \alpha'_L$ in this symmetric Hopf algebroid. Then, from $\beta_R' \circ \epsilon_R' \circ \alpha_L' = \alpha_L'$, we have $\beta'_R \circ \phi = \alpha'_L$, that is
		\begin{equation} \label{eqa}
		\phi(x)_{[-1]} \sharp \phi(x)_{[0]} = \beta'_R(\phi(x)) = \alpha'_L (x)= \theta^{-1}(x_{[0]}) \cdot S(x_{[1]}). 
		\end{equation}
	From this, it follows that the map $\Psi (x \sharp f) := \phi(x)_{[-1]} \sharp \phi(x)_{[0]} \cdot f$ is equal to the inverse of $\bar\Psi$, the map $\bar\Psi^{-1}(x\sharp f) = \theta^{-1}(x_{[0]}) \cdot S(x_{[1]}) f$, where $x\in L$ and $f\in H$. We conclude that $\Psi$ is a homomorphism, thus we have
		$$f \cdot \Psi(x) = \Psi(f \cdot x) = \Psi((f_{(1)} \btr x) f_{(2)}) = \Psi(f_{(1)} \btr x) \cdot f_{(2)}. $$
	We apply $\epsilon'_R$ to the equivalent equation $\Psi(f\btr x) = f_{(1)} \cdot \Psi(x) \cdot S(f_{(2)})$ and  we get the compatibility of actions via $\phi$,
		\begin{equation}
		\phi(f \btr x) = \phi(x) \btl Sf \label{compac}.
		\end{equation}

	Next, compatibility of coactions via $\phi$ is equivalent to the equality $$\theta^{-1}(x_{[0]}) \cdot S(x_{[1]}) = S(x_{[1]}) \sharp \phi(x_{[0]}), \quad \text{ for all } x\in L,$$ since (\ref{eqa}). Set $y:= \theta^{-1}(x)$. The following equations are all equivalent to it, which follows from the above proven compatibility of actions via $\phi$ and the assumed compatibility of coactions via $\theta$,  
		\begin{align*}
		\theta^{-1}(\theta(y)_{[0]}) \cdot S(\theta(y)_{[1]}) &= S(\theta(y)_{[1]}) \sharp \phi(\theta(y)_{[0]}) & \\
		y_{[0]} \cdot S^2(y_{[-1]}) &= S^2(y_{[-1]}) \sharp \phi(\theta(y_{[0]})) & \text{by (\ref{bbb})}\\
		y &= S^2(y_{[-1]}) \sharp \phi(\theta(y_{[0]}))  \cdot S(y_{[-2]}) &\\
		y &= \phi(\theta(y_{[0]}))  \btl S(y_{[-1]}) &\\
		y &= \phi(y_{[-1]} \btr \theta(y_{[0]}))  & \text{by (\ref{compac})}\\
		y & = \phi (\epsilon_L (\beta_R(y))). &
		\end{align*}
	The last equation is actually $\phi^{-1}(y) = \epsilon_L \circ \beta_R$, which holds in a symmetric Hopf algebroid. Therefore, the compatibility of coactions via $\phi$ also holds. 
	
	The corresponding symmetric Hopf algebroid defined from compatibility $\phi$ by using isomorphism $\Psi = \bar\Psi^{-1}$ from Theorem \ref{HAHA} is the same as the one defined from compatibility $\theta$ by using isomorphism $\bar\Psi$ in Theorem \ref{barPhi} because the scalar extension bialgebroid structure maps are the same and the antipodes are the same.
		
	Part (\ref{71}) is proven analogously.	
				
\end{proof}

\subsection{Isomorphism of smash product algebras. Converse theorem} \label{sectionC}
	
	
Given $H$-module algebra structures $(L,\btr)$ and $(R,\btl)$ and $H$-comodule structures $(L,\rho)$ and $(R,\lambda)$ with coactions that are antihomomorphisms, if the corresponding smash product algebras $L\sharp H$ and $H\sharp R$ are isomorphic, the Yetter--Drinfeld conditions in them automatically hold and compatibilites of actions via certain antiisomorphisms~$\phi$ and~$\theta$ also hold. If furthermore the antipodes agree via this isomorphism, compatibilities of coactions via $\phi$ and via $\theta$ hold. This is the content of the Theorem~\ref{obrat} following the next proposition.

\begin{proposition} \label{obrat-1} 
	Let $H$ be a Hopf algebra with an antipode $S$. 
	Let $(L,\btr)$ be a left $H$-module algebra and $(L,\rho)$ a right $H$-comodule. Let $(R, \btl)$ be a right $H$-module algebra  and $(R,\lambda)$ a left $H$-comodule. 
	Assume that the coactions $\rho \colon L \to L\sharp H$ and  $\lambda \colon R \to H \sharp R$ are antihomomorphisms. 	
	
	Let $\phi \colon L \to R$ and $\theta \colon R \to L$ be antiisomorphisms of algebras.

	Denote by $\Psi$ and $\bar\Psi$ the following maps
		\begin{align*}  
			\Psi\colon L \sharp H \to H\sharp R, \quad \Psi(x\sharp f) = \lambda(\phi(x)) \cdot f ,\\
			\bar\Psi \colon H\sharp R \to L \sharp H, \quad \bar\Psi (f\sharp y) =  f \cdot \rho(\theta(y)). 
		\end{align*} 
		\begin{center}
			\begin{tikzcd}
			L \arrow[d, hook, swap]  \arrow[r, "\phi "] & R \arrow[d, "\lambda"]\\
			L\sharp H \arrow[r, dashed, shift left, "\Psi"] & H\sharp R  \arrow[l, dashed, shift left , "\overline \Psi"]\\
			L \arrow[u, "\rho" ] & R \arrow[l, "\theta"] \arrow[u, hook] 
			\end{tikzcd}
		\end{center} 
	Define $\tau \colon L \sharp H \to L \sharp H$ and $\tau' \colon H\sharp R \to H\sharp R$ by 
		\begin{align*}
		\tau ( x \sharp f) & = Sf S^2(x_{[1]}) \cdot x_{[0]}, \\
		\tau' ( f \sharp y) & =   y_{[0]} \cdot S^2(y_{[-1]})Sf.
		\end{align*} 
	Then the following holds.
	\begin{enumerate}
		\item If $\Psi$ is a homomorphism of algebras, then actions are compatible via $\phi$, \label{81}
			$$\phi(f \btr x) = \phi(x) \btl Sf, \quad \text{ for all } x\in L, f\in H.$$ 
		If $\bar\Psi$ is a homomorphism of algebras, then  actions are compatible via $\theta$,
			$$\theta(y\btl f) = Sf \btr \theta(y), \quad \text{ for all } y\in R, f\in H.$$ 
		\item If $\Psi$ is an isomorphism of algebras and $\bar\Psi$ is a homomorphism, then the left-right Yetter--Drinfeld condition holds, that is, for all $x\in L$ and $f \in H$, \label{82}
			$$f \cdot x_{[0]}  \sharp x_{[1]} = (f_{(2)} \btr x)_{[0]}  \sharp (f_{(2)} \btr x)_{[1]}  f_{(1)}.$$
		If $\bar\Psi$ is an isomorphism of algebras and $\Psi$ is a homomorphism, then the right-left Yetter--Drinfeld condition holds,  that is, for all $y\in R$ and $f \in H$,
			$$y_{[-1]} \sharp y_{[0]} \cdot f = f_{(2)}(y \btl f_{(1)})_{[-1]} \sharp (y \btl f_{(1)})_{[0]}.$$
		\item \label{83} If $\Psi$ is an isomorphism and $\tau$ and $\tau'$ agree via $\Psi$, that is $\Psi \circ \tau = \tau' \circ \Psi$, then coactions are compatible via $\phi$,
			\begin{align*}
			\phi(x)_{[-1]}\otimes \phi(x)_{[0]} = S(x_{[1]}) \otimes \phi(x_{[0]}), \quad \text{ for all } x\in L.
			\end{align*}			
		If $\bar\Psi$ is an isomorphism and $\tau$ and $\tau'$ agree via $\bar\Psi$, that is $\bar \Psi \circ \tau' = \tau \circ \bar\Psi$, then coactions are compatible via $\theta$,
			\begin{align*}
			\theta(y)_{[0]}\otimes \theta(y)_{[1]} = \theta(y_{[0]}) \otimes S(y_{[-1]}), \quad \text{ for all } y\in R.
			\end{align*}
		\item \label{84} If $\Psi$ and $\bar \Psi$ are isomorphisms inverse to each other and $\tau$ and $\tau'$ agree, then the symmetric Hopf algebroid defined via $\phi$ and the symmetric Hopf algebroid defined via $\theta$ are the same. 
			
	\end{enumerate}
\end{proposition}

\begin{proof}
	We prove (\ref{81}). If $\Psi$ is a homomorphism, then for all $x\in L$ and $f\in H$ $$\Psi(f\cdot x) = \Psi((f_{(1)} \btr x) \sharp f_{(2)}) = \Psi(f_{(1)} \btr x) \cdot f_{(2)}$$ equals $f \cdot \Psi(x)$. This is equivalent to 
		$\Psi(f \btr x) = f_{(1)} \cdot \Psi(x) \cdot S(f_{(2)}).$
	We apply $\epsilon'_R$ to it and get the compatibility of actions via $\phi$. Similarly, if $\bar \Psi$ is a homomorphism, then for all $y\in R$ and $f\in H$ we have
		$$\bar \Psi (y \btl f) = S(f_{(1)}) \cdot \Psi(y) \cdot f_{(2)}.$$ 
	By applying $\epsilon_L$ to it, we get the compatibility of actions via $\theta$.
		
	We prove (\ref{82}). If $\Psi$ is an isomorphism and $\bar\Psi$ is a homomorphism, we  apply $\Psi$ to the left-right Yetter--Drinfeld condition in $L\sharp H$ $$f \cdot x_{[0]}  \sharp x_{[1]} = (f_{(2)} \btr x)_{[0]}  \sharp (f_{(2)} \btr x)_{[1]}  f_{(1)}$$ and use the compatibility $\theta(y \btl f) = Sf \btr \theta(y)$, which follows from $\bar\Psi$ being a homomorphism, to get a trivially true statement in $H\sharp R$ equivalent to it,
		\begin{align*}
			f \cdot \Psi(\rho(x)) &= \Psi(\rho(f_{(2)} \btr x)) \cdot f_{(1)} \\ f \cdot \Psi (\bar\Psi (1\sharp \theta^{-1}(x))) &=  \Psi (\bar\Psi (1\sharp \theta^{-1}(f_{(2)} \btr x))) \cdot f_{(1)} \\
			(\Psi \circ \bar\Psi) (f \sharp \theta^{-1}(x))) & = (\Psi \circ \bar\Psi) (f_{(1)} \sharp (\theta^{-1}(f_{(3)} \btr x) \btl f_{(2)}))) \\
			(\Psi \circ \bar\Psi)(f \sharp \theta^{-1}(x)) & = (\Psi \circ \bar\Psi)(f_{(1)} \sharp \theta^{-1}(S(f_{(2)}) \btr (f_{(3)} \btr x))) \\
			(\Psi \circ \bar\Psi)(f \sharp \theta^{-1}(x)) & = (\Psi \circ \bar\Psi)(f \sharp \theta^{-1}(x)).
		\end{align*} 
	Similarly, if $\bar\Psi$ is an isomorphism and $\Psi$ is a homomorphism, we  apply $\bar\Psi$ to the right-left Yetter--Drinfeld condition in $H\sharp R$
		$$y_{[-1]} \sharp y_{[0]} \cdot f = f_{(2)}(y \btl f_{(1)})_{[-1]} \sharp (y \btl f_{(1)})_{[0]}$$
	and use the compatibility $\phi(f\btr x) = \phi(x) \btl Sf$, which follows from $\Psi$ being a homomorphism, to get an equivalent equation that is trivially true in $L\sharp H$,
		\begin{align*}
			\bar\Psi(\lambda(y)) \cdot f &= f_{(2)} \cdot \bar\Psi(\lambda(y \btl f_{(1)})) \\
			\bar\Psi (\Psi (\phi^{-1}(y)\sharp 1)) \cdot f &=  f_{(2)} \cdot \bar\Psi (\Psi ( \phi^{-1}(y \btl f_{(1)}) \sharp 1)) \\
			(\bar\Psi \circ \Psi )(\phi^{-1}(y) \sharp f) &= (\bar\Psi \circ \Psi ) ((f_{(2)} \btr \phi^{-1}(y \btl f_{(1)})) \sharp f_{(3)} )\\
			(\bar\Psi \circ \Psi )(\phi^{-1}(y) \sharp f) &= (\bar\Psi \circ \Psi )(\phi^{-1}((y \btl f_{(1)}) \btl S(f_{(2)})) \sharp f_{(3)} )\\
			(\bar\Psi \circ \Psi )(\phi^{-1}(y) \sharp f) &= (\bar\Psi \circ \Psi )(\phi^{-1}(y) \sharp f).
		\end{align*}
		
	We prove (\ref{83}). If $\Psi$ is an isomorphism, then 
		\begin{align*}
		(\Psi \circ \tau)(x\sharp 1) &= \Psi(S^2(x_{[1]}) \cdot x_{[0]}) = S^2(x_{[1]}) \cdot  \Psi(x_{[0]})\\
		&= S^2(x_{[1]}) \cdot \phi(x_{[0]})_{[-1]} \sharp \phi(x_{[0]})_{[0]}, \\
		(\tau' \circ \Psi )(x\sharp 1) &= \tau'(\lambda(\phi(x))) \\
		&= \phi(x).
		\end{align*}
	The equation $(\Psi \circ \tau)(x\sharp 1) = (\tau' \circ \Psi )(x\sharp 1)$ is  therefore 
		$$S^2(x_{[1]}) \cdot \phi(x_{[0]})_{[-1]} \sharp \phi(x_{[0]})_{[0]} = \phi(x),$$which is equivalent to  the compatibility of coactions via $\phi$. 
				
	The equation $(\bar\Psi  \circ  \tau')(1\sharp y)=(\tau \circ \Psi)(1\sharp y)$ is similarly equivalent to the compatibility of coactions via $\theta$.
		
	We prove (\ref{84}). Since both $\Psi $ and $\bar\Psi$ are isomorphisms, both Yetter--Drinfeld conditions are satisfied by (ii), hence we have two scalar extension bialgebroids. By (\ref{83}), if $\tau$ and $\tau'$ agree as in $\Psi \circ \tau = \tau' \circ \Psi$, or as in  $\bar \Psi \circ \tau' = \tau \circ \bar\Psi$, the compatibility of coactions holds via $\phi$ or via $\theta$, respectively, and we have the corresponding symmetric Hopf algebroid structure via $\phi$ or via $\theta$, respectively. These symmetric Hopf algebroid structures are the same if $\Psi$ is inverse to $\bar\Psi$, since the scalar extension bialgebroid structures agree.
\end{proof}

Now we formulate the converse theorem. 
	
\begin{theorem}\label{obrat}
	Let $H$ be a Hopf algebra. Let $(L,\btr)$ be a left $H$-module algebra and $(L,\rho)$ a right $H$-comodule. Let $(R, \btl)$ be a right $H$-module algebra  and $(R,\lambda)$ a left $H$-comodule. Assume that the coactions $\rho \colon L \to L\sharp H$ and $\lambda \colon R \to H \sharp R$ are antihomomorphisms. 
		
	Let $\Psi \colon L\sharp H \to H\sharp R $ be an isomorphism of algebras such that $\Psi(f\sharp 1) = 1\sharp f$ for all $f\in H$, $\Psi(L\sharp 1) = \operatorname{Im} \lambda$ and $\Psi^{-1}(1\sharp R) = \operatorname{Im} \rho$. 
		
	Denote by $\phi \colon L \to R$ and $\theta  \colon R \to L$ antiisomorphisms such that the next diagram commutes.
		\begin{center}
			\begin{tikzcd}
			L \arrow[d, hook]  \arrow[r, dashed, "\phi "] & R \arrow[d, "\lambda"]\\
			L\sharp H \arrow[r, shift left, "\Psi"] & H\sharp R  \arrow[l, shift left , "\Psi^{-1}"]\\
			L \arrow[u, "\rho" ] & R \arrow[l, dashed, "\theta"] \arrow[u, hook] 
			\end{tikzcd}
		\end{center}
		
	Then Yetter--Drinfeld conditions for $(L, \btr, \rho)$ and $(R, \btl, \lambda)$ hold and the actions are compatible via $\phi$ and via $\theta$. 
		
	If furthermore antipodes $\tau$ and $\tau'$ for the corresponding scalar extension bialgebroids $L\sharp H$ and $H\sharp R$ agree, that is $\Psi \circ \tau = \tau' \circ \Psi$, then the coactions are also compatible via $\phi$ and via $\theta$. The resulting symmetric Hopf algebroid structure defined by formulas via $\phi$ from Theorem \ref{HAHA} and the one defined by formulas via $\theta$ from Theorem \ref{barPhi} are the same.
 \end{theorem}
	
\begin{proof}
	Together the braided commutativity and comodule algebra property for $L$ and $R$ follow from the antihomomorphism property of $\rho$ and $\lambda$ respectively. Since $\Psi$ is an isomorphism of algebras, by using definitions of $\phi$ and $\theta$, we have $\Psi(x\sharp f) = \lambda(\phi(x)) \cdot f$ and $\Psi^{-1}(f\sharp y) = f \cdot \rho(\theta(y))$, for all $x\in L$, $f\in H$ and $y\in R$. By Proposition \ref{obrat-1}, actions $\btl$ and $\btr$ are compatible via $\phi$ and via $\theta$ and the Yetter--Drinfeld properties for $(L,\btr,\rho)$ and $(R,\btl,\lambda)$ hold. Therefore, we have scalar extension structures of a left bialgebroid $(L\sharp H,\alpha_L,\beta_L,\Delta_L,\epsilon_L)$ and a right bialgebroid $(H\sharp R,\alpha'_R,\beta'_R,\Delta'_R,\epsilon'_R)$ by Proposition~\ref{bial}. 
		
	Now assume that the antipodes agree. By Proposition \ref{obrat-1}, with $\bar\Psi := \Psi^{-1}$, it follows that the coactions are compatible via $\phi$ and via $\theta$, and two symmetric Hopf algebroid structures arise as a result by Theorems \ref{HAHA} and~\ref{barPhi}. By Proposition \ref{obrat-1}, these structures are the same because $\bar\Psi$ is the inverse of~$\Psi$. 
\end{proof}

\subsection{When the Hopf algebra antipode is bijective} \label{sectionD}

When the Hopf algebra $H$ antipode is bijective, a left-right braided commutative Yetter--Drinfeld $H$-module algebra automatically gives rise to a right-left one compatible to it, via $\phi$ or via $\theta$, and the other way round.	

\begin{proposition}\label{paired} 
	Let $H$ be a Hopf algebra with a bijective antipode $S$. 
		
	Let $\phi \colon L \to R$ be an antiisomorphism of algebras. 
		\begin{enumerate}
			
			\item \label{91} Let $(L, \btr, \rho)$ be a left-right braided commutative Yetter--Drinfeld $H$-module algebra. Then $(R,\btl, \lambda)$ with structure maps defined by
				\begin{align*}
				y \btl f &:= \phi(S^{-1}f \btr \phi^{-1}(y)),
				\\
				y_{[-1]} \otimes y_{[0]} &:= S(\phi^{-1}(y)_{[1]}) \otimes \phi(\phi^{-1}( y)_{[0]} ) , \quad \text{ for } y\in R, f\in H
				\end{align*} 
			is a right-left braided commutative Yetter--Drinfeld $H$-module algebra. 
			
			\item \label{92} Let $(R, \btl, \lambda)$ be a right-left braided commutative Yetter--Drinfeld $H$-module algebra. Then $(L,\btr, \rho)$ with structure maps defined by
			\begin{align*}
			f \btr x &:= \phi^{-1}(\phi(x) \btl Sf) ,
			\\
			x_{[0]} \otimes x_{[1]} &:=  \phi^{-1}(\phi(x)_{[0]}) \otimes  S^{-1}(\phi(x)_{[-1]}), \quad \text{ for } x\in L, f\in H
			\end{align*} 
			is a left-right braided commutative Yetter--Drinfeld $H$-module algebra.

		\end{enumerate}	
	
	Let $\theta \colon R \to L$ be an antiisomorphism of algebras.
		\begin{enumerate} \setcounter{enumi}{2}

			\item \label{93} Let $(L, \btr, \rho)$ be a left-right braided commutative Yetter--Drinfeld $H$-module algebra. Then $(R,\btl, \lambda)$ with structure maps defined by
				\begin{align*}
				y \btl f &:= \theta^{-1}(S f \btr \theta(y)),
				\\
				y_{[-1]} \otimes y_{[0]} &:= S^{-1}(\theta(y)_{[1]}) \otimes \theta^{-1}( \theta( y)_{[0]})  , \quad \text{ for } y\in R, f\in H
				\end{align*} 
			is a right-left braided commutative Yetter--Drinfeld $H$-module algebra. 
			
			\item \label{94} Let $(R, \btl, \lambda)$ be a right-left braided commutative Yetter--Drinfeld $H$-module algebra. Then $(L,\btr, \rho)$ with structure maps defined by
				\begin{align*}
				f \btr x &:= \theta(\theta^{-1}(x) \btl S^{-1}f) 
				\\
				x_{[0]} \otimes x_{[1]} &:= \theta( \theta^{-1}(x)_{[0]}) \otimes  S(\theta^{-1}(x)_{[-1]}), \quad \text{ for } x\in L, f\in H
				\end{align*} 
			is a left-right braided commutative Yetter--Drinfeld $H$-module algebra.

		\end{enumerate}	
		
\end{proposition}
	
\begin{proof}
	We prove (\ref{91}) and (\ref{92}). It is easy to check that  $\btr$ and $\lambda$ are Hopf action and coaction if and only if $\btl$ and $\rho$ are, respectively. Furthermore, it is easy to see that right-left Yetter--Drinfeld condition for $Sf$ and $\phi(x)$ is equivalent to $\phi\otimes S$ of a left-right Yetter--Drinfeld condition for $f\in H$ and $x \in L$. From this, since $S$ and $\phi$ are bijective, equivalence of two Yetter--Drinfeld conditions follows. The corresponding  braided commutativity properties and comodule algebra properties are easily proven to be equivalent by using antiisomorphism $\phi$. The proof of (\ref{93}) and (\ref{94}) is analogous. 	
\end{proof}	
	
\begin{corollary} 
	Let $H$ be a Hopf algebra with a bijective antipode $S$. Let $\phi \colon L \to R$ be an antiisomophism of algebras. 
		
	\begin{enumerate} 
	
		\item \label{101} Let $(L,\btr, \rho)$ be a left-right braided commutative Yetter--Drinfeld $H$-module algebra. Then $L\sharp H$ with the scalar extension structure maps for the left bialgebroid, the following structure maps for the right bialgebroid and antipode as below is a symmetric Hopf algebroid. 
			\begin{align*} 
			\alpha_R  \colon& R \to L\sharp H &&& \alpha_R( y ) &=  S^2(\phi^{-1}(y)_{[1]}) \cdot \phi^{-1}(y)_{[0]}\\
			\beta_R  \colon & R \to L\sharp H &&& \beta_R( y ) &= \phi^{-1}(y) \sharp 1_H \\ 
			\Delta_R  \colon & L \sharp H \to L\sharp H \otimes_{R} L \sharp H &&& \Delta_R( x \sharp f ) &= x \sharp f_{(1)} \otimes_R 1_L \sharp f_{(2)}\\
			\epsilon_R \colon& L \sharp H \to R &&& \epsilon_R(x \sharp f) &=  \phi( S^{-1}f \btr x) \\
			\tau  \colon & L \sharp H \to L \sharp H &&& \tau(x \sharp f) &= Sf S^2(x_{[1]}) \cdot x_{[0]} .
			\end{align*}
			
		\item \label{102} Let $(R,\btl,\lambda)$ be a right-left braided commutative Yetter--Drinfeld $H$-module algebra. Then $H\sharp R$ with the scalar extension structure maps for the right bialgebroid, the  following structure maps for the left bialgebroid and antipode as below is a symmetric Hopf algebroid.
			\begin{align*} 
			\alpha_L' \colon& L \to H\sharp R &&& \alpha_L'(x) &= \lambda(\phi(x)) = \phi(x)_{[-1]} \sharp  \phi(x)_{[0]} \\
			\beta_L'  \colon& L \to H\sharp R &&& \beta_L'(x) &= 
			\phi(x)_{[0]} \btl S^{-1} (\phi(x)_{[-1]})\\ 
			\Delta_L'  \colon & H \sharp R \to H\sharp R \otimes_{L} H \sharp R &&& \Delta_L'(f \sharp y) &= f_{(1)}\sharp 1_R \otimes_L f_{(2)}\sharp y \\
			\epsilon_L' \colon& H \sharp R \to L &&& \epsilon_L'(f\sharp y) &= \phi^{-1}(y_{[0]}\btl S^2(y_{[-1]})Sf )
			\\
			\tau'  \colon & H \sharp R \to H \sharp R &&& \tau'(f \sharp y ) & = y_{[0]} \cdot S^2(y_{[-1]})Sf.
			\end{align*}
	\end{enumerate}
	Furthermore, if one Yetter--Drinfeld algebra is determined by the other as in Proposition \ref{paired}, these are two presentations of the same symmetric Hopf algebroid, connec\-ted by isomorphism $\Psi \colon L\sharp H \to H\sharp R$ defined in Theorem~\ref{Phi}. 
\end{corollary}

Analogous proposition with antiisomorphism $\theta$ is also a corollary of the theorems in the previous sections. All formulas for symmetric Hopf algebroid structure maps, when antipode $S$ is bijective, are organized in the following table.

	\begin{align*}
	& \text{From datum $(L,\btr,\rho)$:}  &  & \text{From datum $(R,\btl,\lambda)$:} \\ \\
	& L\sharp H&        &             H \sharp R  \\
	& \emph{left bialgebroid over } L  & & \text{left bialgebroid over }  L, \text{ by } \phi \\ 
	& \alpha_L(x)= x \sharp 1_H &              &                 \alpha_L(x) =  \phi(x)_{[-1]} \sharp \phi(x)_{[0]}  \\
	&\beta_L(x) = x_{[0]} \sharp x_{[1]}          & & \beta_L(x) =    \phi(x)_{[0]} \btl S^{-1} (\phi(x)_{[-1]}) \\
	& \Delta_L(x \sharp f) = x \sharp f_{(1)} \otimes_L 1_L \sharp f_{(2)}      &  & \Delta_L (f \sharp y) =  f_{(1)} \sharp 1_R \otimes_L f_{(2)} \sharp y  \\
	& \epsilon_L(x \sharp f) =\epsilon(f) x &   &   \epsilon_L(f\sharp y)  =  \phi^{-1} (y_{[0]} \btl  S^{2}(y_{[-1]}) Sf)  \\
	\\
	& L\sharp H&        &             H \sharp R  \\
	& \text{right bialgebroid over } R, \text{ by } \phi & & \emph{right bialgebroid over }  R  \\   
	& \alpha_R(y)= S^2(\phi^{-1}(y)_{[1]}) \cdot \phi^{-1}(y)_{[0]} & &                              \alpha_R(y) =  1_H \sharp y   \\
	& \beta_R(y) = \phi^{-1}(y) \sharp 1_H          & & \beta_R(y) = y_{[-1]} \sharp y_{[0]}  \\
	&  \Delta_R(x \sharp f) = x \sharp f_{(1)} \otimes_R 1_L \sharp f_{(2)}      &   & \Delta_R(f \sharp y) = f_{(1)} \sharp 1_R \otimes_R f_{(2)} \sharp y   \\
	& \epsilon_R(x \sharp f) = \phi(S^{-1}f  \btr x)  & &   \epsilon_R(f \sharp y) = \epsilon(f)y   \\
	\\
	& \text{antipode}  &  & \text{antipode}  \\ 
	& \tau (x \sharp f) =  Sf S^2(x_{[1]}) \cdot x_{[0]} & &  \tau (f \sharp y) = y_{[0]} \cdot S^2(y_{[-1]})Sf  \\
	& \tau^{-1} (x \sharp f) = S^{-1}f \cdot  x_{[0]} \sharp x_{[1]}  & & \tau^{-1} (f \sharp y) = y_{[-1]} \sharp y_{[0]}  \cdot S^{-1}f \\
	\\
	& \text{Here }\phi = \epsilon_R \circ \alpha_L,\ \phi^{-1} = \epsilon_L \circ \beta_R,\ & & \phi \colon L \to R \text{ antiisomorphism.} 
	\end{align*}
	Alternatively, formulas  can be written by using the antiisomorphism $\theta$. 
	\begin{align*}
	& L\sharp H&                    & H \sharp R\\
	&\text{right bialgebroid over } R, \text{ by } \theta & & \text{left bialgebroid over }  L, \text{ by } \theta \\   
	&\alpha_R(y)= \theta(y)_{[0]} \sharp \theta(y)_{[1]} &   &       \alpha_L(x) =   \theta^{-1}(x)_{[0]} \cdot S^2(\theta^{-1}(x)_{[-1]})   \\
	&\beta_R(y) =   S^{-1} (\theta(y)_{[1]}) \btr \theta(y)_{[0]}  & & \beta_L(x) =    1_H \sharp \theta^{-1}(x)     \\
	& \Delta_R(x \sharp f) = x \sharp f_{(1)} \otimes_R 1_L \sharp f_{(2)}      &  & \Delta_L(f \sharp y) = f_{(1)} \sharp 1_R \otimes_L f_{(2)} \sharp y   \\
	&\epsilon_R(x \sharp f) =  \theta^{-1} ( SfS^2(x_{[1]}) \btr x_{[0]}) & & \epsilon_L(f\sharp y)  =  \theta( y \btl S^{-1}f )  \\
	\\
	& \text{Here } \theta = \epsilon_L \circ \alpha_R, \ \theta^{-1} = \epsilon_R \circ \beta_L, & & \theta \colon  R \to L \text{ antiisomorphism.} 
	\end{align*}

\section*{Acknowledgements}  The author is grateful to Zoran Škoda for discussions, encouragement and for helping with the literature.

\end{document}